\documentclass{article}

\PassOptionsToPackage{numbers, compress}{natbib}

\usepackage{definitions}

\usepackage[final]{neurips_2023}

\usepackage[utf8]{inputenc} 
\usepackage[T1]{fontenc}    
\usepackage{hyperref}       
\usepackage{url}            
\usepackage{booktabs}       
\usepackage{amsfonts}       
\usepackage{nicefrac}       
\usepackage{microtype}      
\usepackage{xcolor}         
\usepackage{graphicx}
\usepackage{comment}
\usepackage{mathrsfs}
\usepackage{amssymb}
\usepackage{amsthm}
\usepackage{amsmath}
\usepackage{amsfonts}
\newcommand{\mathbbm}[1]{\text{\usefont{U}{bbm}{m}{n}#1}}
\newtheorem{theorem}{Theorem}

\newtheorem{assumption}{Assumption}
\newtheorem{lemma}{Lemma}
\theoremstyle{remark}
\newtheorem{remark}{Remark}

\title{Optimal testing using combined test statistics across independent studies}

\author{%
  Lasse~Vuursteen \\
  Delft Institute of Applied Mathematics\\
  Delft University of Technology\\
  \texttt{l.vuursteen@tudelft.nl} \\
  \And
  Botond~Szab\'o \\
  Department of Decision Sciences \\
  and Institute for Data Science and Analytics \\
  Bocconi University \\
  \texttt{botond.szabo@unibocconi.it} \\
   \AND
  Aad van der Vaart \\
  Delft Institute of Applied Mathematics\\
  Delft University of Technology\\
  \texttt{a.w.vandervaart@tudelft.nl} \\
  \And
  Harry van Zanten \\
  Mathematics Department \\
  Vrije Universiteit Amsterdam \\
  \texttt{j.h.van.zanten@vu.nl} \\
}

\begin{document}

\maketitle

\begin{abstract}
Combining test statistics from independent trials or experiments is a popular method of meta-analysis. However, there is very limited theoretical understanding of the power of the combined test, especially in high-dimensional models considering composite hypotheses tests. We derive a mathematical framework to study standard {meta-analysis} testing approaches in the context of the many normal means model, which serves as the platform to investigate more complex models.

We introduce a natural and mild restriction on the meta-level combination functions of the local trials. This allows us to mathematically quantify the cost of compressing $m$ trials into real-valued test statistics and combining these. We then derive minimax lower and matching upper bounds for the separation rates of standard combination methods for e.g. p-values and e-values, quantifying the loss relative to using the full, pooled data. We observe an elbow effect, revealing that in certain cases combining the locally optimal tests in each trial results in a sub-optimal {meta-analysis} method and develop approaches to achieve the global optima. We also explore the possible gains of allowing limited coordination between the trial designs. Our results connect meta-analysis with bandwidth constraint distributed inference and build on recent information theoretic developments in the latter field. 
\end{abstract}

\section{Introduction}

Given multiple data sets relating to the same hypothesis, one would like to combine the evidence. Sometimes, the full data sets are not available (e.g. due to privacy or proprietary reasons) or difficult to combine directly (e.g. due to the different experimental or observational setups). In such cases, the analysis must be carried out on the basis of the published results for each of the studies. Such {``meta-analysis''} can increase the statistical power by combining individually inconclusive or moderately significant tests, while keeping the false positive rate under control. Therefore, meta-analysis has received a lot of attention in various fields, for instance in genetics and system biology, when studying rare variants \cite{aerts2006gene,evangelou2013meta} or in deep learning, for few shot image recognition and neural architecture search, see the review article \cite{hospedales2021meta}.

The outcomes of the studies concerning hypothesis tests are, typically, summarized as real-valued test statistics and/or associated p-values. One expects the combination of $m$ such p-values to result in an increase in power, but one also expects to pay a price relative to computing a test on the basis of the full, pooled data of the $m$ trials. The question of how to optimally combine independent real-valued test statistics concerning the same hypothesis into a single test has an extensive literature. A multitude of methods for combining independent tests of significance exist. For combining p-values, this starts with Fisher, Tippett and Pearson in the nineteen-thirties, see \cite{tippett1941methods,fisher1992statistical, pearson1934new,  stouffer1949american, good1955weighted, liptak1958combination, van1967combination, edgington1972additive, mudholkar1977logit, whitlock2005combining, vovk2020combining, chen2021optimal, vovk2022admissible} and references therein. In Section \ref{sec:examples_worked_out}, we collect and describe the most popular and frequently used p-value combination techniques. 

As noted in \cite{birnbaum1954combining}, there does not exist a general uniformly most powerful p-value combination method for all alternative hypotheses. The distribution of a p-value or its underlying test statistic under the alternative hypothesis should be taken into consideration when selecting a method of combination. The performance of different p-value combination techniques was investigated extensively by empirical experiments in various synthetic and real world scenarios, see for instance \cite{loughin2004systematic,yoon2021powerful}. However, a unified, general theoretical description is lacking, especially in non-trivial, multi-dimensional composite testing problems, where the likelihood ratio test is not necessarily uniformly most powerful. 

E-values are an increasingly popular and important notion of evidence, see \cite{shafer_test_2011, IEEE_Safe_Testing,shafer2021testing}. E-values allow the combination of several tests in a straightforward manner while preserving the prescribed level of the tests (see Section \ref{ssec:e-values_combined}). Formally, e-values are nonnegative random variables whose expected values under the null hypothesis are bounded by one. In contrast to p-values defined by probabilities, e-values are defined by expectation. This imposes significant differences in their interpretation, application and combination compared to the more standard p-values. However, as for p-values, very little is known about the power of these combination procedures. Theoretical results focus on specific optimality criteria, for instance the worst-case growth-rate (GROW), see \cite{IEEE_Safe_Testing}. However, these do not directly imply guarantees on the testing power, which is the main focus in practice.

Our focus is on multidimensional models, where a certain loss in power is to be expected, since combining multidimensional data into a real-valued statistic (e.g. p-value or e-value) requires data compression. Typically, summary statistics are combined by some ``reasonable'' function $C_m:\R^m \to \R$, where ``reasonable'' means that $C_m$ should not exploit the richness of the real numbers to encode the data in full. We aim to quantify the loss of summarizing, the gain of performing a meta-analysis and the best testing strategies in the individual experiments meta-analysis.

We consider the signal detection problem in the many normal means model, see Section \ref{sec:main_results} for the detailed description. One possible interpretation of this testing problem is to learn whether a treatment has an effect on any of the dimensions investigated. {This model is directly applied in several fields where high-dimensional statistics and machine learning settings are concerned, such as detecting differentially expressed genes \cite{quackenbush2002microarray,kramer2014causal,malone2011comparison,thomas2001efficient,efron2012large}, bankruptcy prediction for publicly traded companies using Altman’s Z-score in finance \cite{altman1968financial,altman2014distressed}, separation of the background and source in astronomical images \cite{clements2012astronomical,guglielmetti2009background}, and wavelet analysis \cite{abramovich2006adapting,johnstone2004needles}.} Furthermore, the model allows for tractable computations and it typically serves as the platform to investigate more difficult statistical and learning problems, including high- and infinite-dimensional models, see for instance \cite{ingster_minimax_1987,tsybakov_introduction_2009,efron2012large,gine_mathematical_2016}. In each experiment $j\in\{1,...,m\}$ the observations are summarized by an appropriate real-valued summary statistic $S^{(j)}$. These local test statistics (e.g. p- or e-values) are {combined} into $C_m(S^{(1)},\dots,S^{(m)})$. We consider a general class of combination functions $C_m$, requiring only H\"older type continuity. This introduces only a mild restriction, and includes many standard {meta-analysis} techniques, for instance the standard p-value combination methods (see Section \ref{ssec:p-values_combined}); e-value techniques (see Section \ref{ssec:e-values_combined}); and other ad hoc and natural test statistic combination approaches, see the beginning of Section \ref{sec:examples_worked_out} for additional examples.
 
Our setting provides a principled and unified framework to study the power of standard meta-analysis testing methods. {Within the framework of the many normal means model,} we derive a minimax lower bound for the testing (separation) error and provide test statistics with associated combination methods that attain this theoretical limit (up to a logarithmic factor). Our results reveal that there is a certain unavoidable loss associated with compressing the data of each experiment to a real valued test statistic. We see that while it is always possible to obtain better testing rates using $m$ trials instead of the best possible test based on a single trial, there is always a loss incurred when compared to the full, pooled data and optimal test in moderate- to large dimensional problems. Our theoretical results quantify these gains and losses in terms of the dimension $d$, sample size $n$ and number of trials $m$. 

Furthermore, we observe an elbow effect, which occurs when the number of trials is large compared to the dimension of the signal. {In this regime, combinations of the (locally) optimal test in each individual trial performs sub-optimally as a whole when aggregated and meta-analysis approaches based on directional test statistics are shown to perform better.} Finally, we show that the performance of the meta-level tests can substantially improve (in certain regimes, depending on $d,m,n$) if a certain amount of coordination between the trials is allowed (e.g. by having access to the same random seed). For the theoretical analysis of {meta-analysis} techniques we derive connections with the distributed statistical learning literature under communication constraints. Our paper builds on the recent information theoretical developments in distributed testing \cite{acharya2020distributed,szabo2022optimal_IEEE,szabo2023optimal}, allowing us to address several fundamental questions for the first time with mathematical rigor.

The paper is organised as follows. In Section \ref{sec:main_results} we introduce the mathematical framework we consider in our investigation and present the corresponding minimax testing lower bound results. Next in Subsection \ref{ssec:upperbound_private_coin} we show that the derived results are sharp by providing several {meta-analysis} approaches attaining the limits. Then we investigate the benefits of allowing a mild coordination between the trials in Subsection \ref{ssec:public_coin}. We collect and discuss the standard p- and e-value combination methods in Section \ref{sec:examples_worked_out} and demonstrate our theoretical results numerically on synthetic data sets in Section \ref{sec:sim}. We discuss our results and derive conclusions in Section \ref{sec:discussion}. The proofs of our results are deferred to the Appendix. In Section \ref{ssec:main_theorem_proof} we present the proof of our main results while the proofs of the technical lemmas are given in \ref{ssec:auxilliary_lemmas}.

 \textbf{Notation:} For two positive sequences $a_n$, $b_n$ we write $a_n\lesssim b_n$ if the inequality  $a_n \leq Cb_n$ holds for some universal positive constant $C$. Similarly, we write $a_n\asymp b_n$ if $a_n\lesssim b_n$ and $b_n\lesssim a_n$ hold simultaneously and let $a_n\ll b_n$ denote that $a_n/b_n =
o(1)$. Furthermore, we use the notations $a\vee b$ and $a\wedge b$ for the maximum and minimum, respectively, between $a$ and $b$. Throughout the paper $c$ and $C$ denote global constants whose value may change from one line to another.

\section{Main results}\label{sec:main_results}

In our analysis, we consider the localized version of the many normal means model, tailored to investigating {meta-analysis} techniques. We assume that in each local trial or experiment $j\in\{1,...,m\}$ we observe a $d$-dimensional random variable $X^{(j)}\in\mathbb{R}^d$, subject to
\begin{equation}\label{dynamics}
X^{(j)}  = f +  \frac{1}{\sqrt{n}} Z^{(j)},\qquad  Z^{(j)}\stackrel{iid}{\sim} N(0,I_d),\quad j=1,...,m,
\end{equation}
for some unknown $f\in\mathbb{R}^d$. We denote by $\P_f$ the joint distribution of the observations {and let $\E_f$ be the corresponding expectation.} We note that this framework is equivalent to having $n$ independent $N(f,I_d)$ observations within each local sample.

Our goal is to test the presence or absence of the ``signal component'' $f\in\mathbb{R}^d$. More formally, we consider the simple null hypothesis $H_0: f = 0$ versus composite alternative hypothesis $H_\rho: \|f\|_2 \geq \rho$, for some $\rho > 0$. This corresponds to testing for joint significance of variables, such as the presence of an effect of a treatment on any of the dimensions investigated. The difficulty in distinguishing the hypotheses depends on the effect size, the sample size and the dimension $d$. Here, $\rho$ can be seen as the smallest effect size deemed important.  

For a $\{0,1\}$-valued test $T$, define the testing risk $\cR(T, H_\rho)$ as the sum of the Type I error probability and worst case Type II error probability, i.e. 
\begin{equation}\label{eq : testing risk}
\cR(T, H_\rho) := \P_0 (T = 1)+ \underset{f \in H_\rho}{\sup} \P_f \left( T = 0 \right).
\end{equation}
In the case of a single trial (i.e. $m=1$), this testing problem is known to have minimax separation rate or ``detection boundary'' $\rho^2 \asymp \sqrt{d}/n$. 

This means that if $\rho^2 \gg \sqrt{d}/n$, there exist consistent\footnote{For any asymptotics in $\rho$, $d$ and $n$ such that $\rho^2 \gg \sqrt{d}/n$.} tests $T \equiv T_{d,n}$ in the sense that $\cR(T, H_\rho) \to 0$, whilst no consistent tests exist when $\rho^2 \ll \sqrt{d}/n$. That is, for effect sizes of smaller order than $\sqrt{d}/n$, the null hypothesis cannot be consistently distinguished from the alternative hypothesis. Such a testing rate is attainable through a chi-square test based on $\|\sqrt{n} X^{(1)}\|_2^2$ (see e.g. \cite{baraud2002non}). 

In case of $m$ trials, if the full data {were} pooled (with aggregated sample size $nm$), the minimax separation rate would be $\sqrt{d}/(mn)$. However, pooling the data might not be possible or allowed in practice and often only real-valued test statistics are available that describe the significance in the local problems (e.g. a p- or an e-value). These $m$ test statistics $S^{(j)}$, $j=1,...,m$, then can be combined with some combination function $C_m:\R^m \to \R$, providing the test statistic in the meta-analysis.  We now ask whether the above pooled testing rate is attainable with this {meta-analysis} procedure.

Without any restrictions on the test statistics $S = (S^{(1)},\dots,S^{(m)})$ or the combination function $C_m$, any of the conventional optimal ``full-data'' tests can be reconstructed, since the real numbers and mappings between the real numbers form an overly rich  class. We wish to restrict our analysis to $S$ and $C_m$ that are reasonable in practice and capture (most of) the relevant meta-analysis methods as listed in Section \ref{sec:examples_worked_out}. 
 
Based on each of the local observations $X^{(j)}$, a real-valued test statistic $S^{(j)}$ is computed, where each $S^{(j)}$ is a function of $X^{(j)}$ and possibly a source of randomness $U^{(j)}$ independent of $X:=(X^{(1)},\dots,X^{(m)})$. 
\begin{assumption}\label{assumption:local_randomness_statistic}
For measurable functions $f_j: \R^d \times \R \to \R$ and independent random variables $U^{(1)},\dots,U^{(m)}$ which are independent of the data $X$, the $j$-th test statistic $S^{(j)} = f_j(X^{(j)},U^{(j)})$ satisfies $\E_0 |S^{(j)}| \leq M$, for some $M>0$, $j=1,\dots,m$.
\end{assumption}
\noindent We consider H\"older continuous combination functions $C_m:\R^m \to \R$. Arguably, this is the most important assumption in ruling out bijections between $\R^d$ and $\R$. This ensures that a small change in the underlying local test statistics cannot result in a large change in the combination of test statistics $C_m(S^{(1)},\dots,S^{(m)})$.
\begin{assumption}\label{assumption:lipschitz_type_C_m}
There exist $L,p,q > 0$ such that for all $s,s' \in \R^m$
\begin{equation}\label{eq:C_m_approximation_assumption}
\left|C_m(s) - C_m(s')\right| \leq L \Big(\underset{j=1}{\overset{m}{\sum}} |s_j - s_j'|^p \Big)^q.
\end{equation}
\end{assumption}
\noindent The special case of $p=2$ and $q=1/2$ leads to Lipschitz continuous functions. {Assumption \ref{assumption:local_randomness_statistic} and Assumption \ref{assumption:lipschitz_type_C_m} should be considered in conjunction.} By rescaling and centering test statistics ${S}^{(j)}$, one can typically obtain test statistics satisfying Assumption \ref{assumption:local_randomness_statistic}. Rescaling and centering typically does affect how the test statistics need to be combined, which might ``break'' Assumption \ref{assumption:lipschitz_type_C_m}. 


Finally, following the standard testing approach, we compare the aggregated test statistics $C_m(S^{(1)},\dots,S^{(m)})$ to a threshold value. If the combined test statistics result in a large enough value, the null hypothesis of no effect is rejected. {We note here that two sided tests can be written as one-sided tests through straightforward transformations (e.g. centering and taking absolute value).} More formally, we consider tests $T_\alpha$ of level $\alpha$ satisfying the following assumption.
\begin{assumption}\label{assumption:test_form}
There exists a strictly decreasing function $\alpha \mapsto \kappa_\alpha$ so that 
\begin{equation}\label{eq:combined_p-value_test}
T_\alpha = \mathbbm{1} \left\{ C_m(S^{(1)},\dots,S^{(m)}) \geq \kappa_{\alpha} \right\}
\end{equation}
satisfies $\P_0 T_\alpha \leq \alpha$. 
\end{assumption}
\noindent The map $\alpha \mapsto \kappa_\alpha$ could be taken as the quantile function of $C_m(S^{(1)},\dots,S^{(m)})$ under its null distribution if it is appropriately standardized. If $\E_0 C_m(S^{(1)},\dots,S^{(m)})$ is bounded in $m$, we can choose $\kappa_\alpha$ equal to $1/\alpha$ times the upper bound, in view of Markov's inequality.

Our first main result, Theorem \ref{thm:private_coin_lb} below, establishes a lower bound for tests of the form \eqref{eq:combined_p-value_test} and $C_m$ and $S$ satisfying the above assumptions.  More concretely, under our assumptions, any test $T_{\alpha}$ (of level $\alpha\leq 0.1$) has large Type II-error under alternatives with $\rho^2$ of smaller order than ${(\sqrt{m}\wedge \frac{d}{\log(m)}) \sqrt{d}}/({m n})$. When the number of trials is small compared to the dimension (i.e. $m \log^2(m) \leq d^2$), this means that the separation rate is at least $\sqrt{d}/(\sqrt{m} n)$. Thus even though there is a benefit in terms of separation rate compared to testing based on just a single trial, the gain is at best the square root of what one would gain based on testing on the pooled data. When $m \log^2(m) \geq d^2$, the rate in the lower bound changes to $d \sqrt{d} / (m n \log (m))$, resulting in an elbow effect.

\begin{theorem}[]\label{thm:private_coin_lb}
Let $S^{(1)},\dots,S^{(m)}$, $C_m$ and $T_{\alpha}$ satisfy Assumptions \ref{assumption:local_randomness_statistic}--\ref{assumption:test_form} with $T_\alpha$ of level $\alpha \in (0,0.1]$. Then there exists a constant $c > 0$ depending only on $L$, $p,q$ and $M$, such that if
\begin{equation}\label{eq:rho_condition_private_coin_lb}
\rho^2 \leq c  \frac{(\sqrt{m}\wedge \frac{d}{\log(m)}) \sqrt{d}}{m n},
\end{equation}
 it holds for all $n,m,d\in \N$ that
\begin{equation}
\underset{f \in H_\rho}{\sup} \P_f \left( T_{\alpha} = 0 \right) \geq 3/4.
\end{equation}
\end{theorem}
\begin{remark}
The ranges of values $0<\alpha \leq 0.1$ and $\beta=3/4$ for the Type I and II errors, respectively, are arbitrary. Similar results hold for different choices as well. For instance, one can take arbitrary $\alpha \in (0,1/5]$ and $\beta \in (0,2/3]$, see the proof of the theorem for details. The result implies in particular that consistent testing is not possible for signals of a smaller order than the right hand side of \eqref{eq:rho_condition_private_coin_lb}, where asymptotics can be considered in $n,m$ and $d$ simultaneously.
\end{remark}

In the next section we show that the lower bounds in the theorems above are sharp (up to a logarithmic factor).

\subsection{Rate optimal combination methods}\label{ssec:upperbound_private_coin}

To attain the lower bound rate derived in Theorem \ref{thm:private_coin_lb}, different tests can be considered. The optimal rate displays an elbow effect around $m \asymp d^2$. When the dimension is large compared to the number of trials $m$ (i.e. $m \lesssim d^2$), strategies that combine p-values for the optimal local tests (based on $\| \sqrt{n} X^{(j)}\|_2^2\sim^{H_0}\chi^2_d$), turn out to achieve the optimal rate, as exhibited below. Such a test statistic is invariant to the directionality of $X^{(j)}$ and invariant under the model in the sense that the resulting power for the alternative $\P_f$ or $\P_g$ is the same as long as $\|f\|_2 = \|g\|_2$. 

On the other hand, when the dimension is small compared to the number of trials (i.e. $m \gtrsim d^2$), optimal strategies exhibited below use information on the direction of $X^{(j)}$. In fact, we show in Theorem \ref{thm:sign_info_necessary_for_rate} in the Appendix that if no such information is available (i.e. the events defined by the signs of the $(X^{(j)})_{j=1,...,m}$ vector are not contained in the sigma algebra generated by the test statistics $S$), one cannot obtain a rate better than ${\sqrt{d}}/({\sqrt{m}n})$. This implies that by combining the locally optimal test statistics  $S^{(j)}=\| \sqrt{n} X^{(j)}\|_2^2$ (or their arbitrary functions, e.g. the corresponding local p-values) would result in information loss and hence sub-optimal rates in the meta-analysis. 

Furthermore, it turns out, in accordance with the empirical literature discussed in the introduction, that there does not exist a uniquely best meta-analysis method. In fact, multiple standard meta-analysis techniques provide (up to a logarithmic factor) optimal rates, see below for some standard approaches attaining the lower bounds derived in Theorem \ref{thm:private_coin_lb}. 

First we consider the scenario when the dimension $d$ of the model is large compared to the number of trials $m$, i.e. $m \lesssim d^2$. Locally the optimal test is based on the test statistic $\| \sqrt{n} X^{(j)}\|_2^2\stackrel{H_0}{\sim}\chi^2_d$. A natural way to combine these statistics would be to sum these locally optimal test statistics to obtain
\begin{equation}
T_\alpha = \mathbbm{1} \left\{ \underset{j=1}{\overset{m}{\sum}} \left\| \sqrt{n} X^{(j)}\right\|_2^2 \geq F_{\chi_{dm}^2}^{-1} (1 - \alpha) \right\},\label{test:large:d}
\end{equation}
which has level $\alpha$. Alternatively, one could also apply p-value combination methods, such as Fisher's or Edgington's method based on the p-value $p^{(j)} = 1 - F_{\chi^2_d}(\| \sqrt{n} X^{(j)}\|_2^2)$, see Section \ref{sec:examples_worked_out}. Lemma \ref{lem:sqrtm_rate_attained} in the appendix establishes that these tests are rate optimal. 

Second, consider the case that the number of trials is large compared to the dimension, i.e. $m \gtrsim d^2$. Rate optimal tests can be constructed based on a variation of Edgington's or Stouffer's method, see Section \ref{sec:examples_worked_out} for their descriptions.
Taking a partition of $\{1,\dots,m\} = \cup_{i=1}^{d} \cJ_i$ where $|\cJ_i| \asymp {m}/{d} $ and setting $S^{(j)} = \sqrt{n^{}} X_i^{(j)}$ if $j \in \cJ_i$, the meta-level test
\begin{equation}\label{eq:private_coin_test}
T_\alpha = \mathbbm{1}\left\{ \frac{\sqrt{d}}{m} \underset{i=1}{\overset{d}{\sum}} \left( \underset{j \in \cJ_i}{\overset{}{\sum}} S^{(j)}  \right)^2  \geq d^{-1/2}F^{-1}_{\chi_d^2}(1-\alpha) \right\}
\end{equation}
achieves the lower bounds. The above test is similar to employing Stouffer's method for each of the coordinates and averaging, i.e. computing approximately $m/d$ iid p-values $p^{(j)} = \Phi(\sqrt{n^{}} X_i^{(j)})$ for $j \in \cJ_i$ and applying the inverse Gaussian CDF $\Phi^{-1}(p^{(j)})$. Alternatively, the following variation of Edgington's method,
\begin{equation}\label{eq:edginton_like_test}
T_\alpha = \mathbbm{1}\left\{ \frac{\sqrt{d}}{m} \underset{i=1}{\overset{d}{\sum}} \left( \underset{j \in \cJ_i}{\overset{}{\sum}}  \left(p^{(j)} - \frac{1}{2}\right)\right)^2  \geq \kappa_\alpha \right\},
\end{equation}
is also rate optimal, as proven in Lemma \ref{lem:priv_coin_coordinate_strat_rate} in the appendix. {Essentially, these strategies divide the trials accross the $d$ different directions, and combines the evidence for each of the directions. Theorem \ref{thm:sign_info_necessary_for_rate} affirms that the information on the ``direction'' of the data is crucial to achieve the optimal rate in the $m \gtrsim d^2$ case, by showing that strategies that do not contain such information (rotationally invariant strategies such as norm-based test statistics) achieve the rate $\sqrt{d}/(\sqrt{m}n)$ at best.} We summarize the above testing upper bounds in the theorem below. 

{\begin{theorem}\label{thm:private_coin_ub}
For all $\alpha,\beta \in (0,1)$ there exist $S$, $C_m:\R^m \to \R$ and tests $T_{\alpha}$ of level $\alpha$ satisfying Assumptions \ref{assumption:local_randomness_statistic}--\ref{assumption:test_form} such that if
\begin{equation}
\rho^2 \geq C_{\alpha,\beta}  \frac{(\sqrt{m}\wedge d) \sqrt{d}}{m n},
\end{equation}
we have
\begin{equation*}
\underset{f \in H_\rho}{\sup} \P_f \left( T_\alpha = 0 \right) \leq \beta
\end{equation*}
for a large enough constant $C_{\alpha,\beta} > 0$ depending only on $\alpha,\beta \in (0,1)$, for all $n,m,d \in \N$.
\end{theorem}}

\subsection{Benefits of coordination between the trials}\label{ssec:public_coin}

When the dimension is small relative to the number of trials, as exhibited in the previous section, optimal strategies include information on the directionality of the observation vector. In this section we show that in this regime, there could be an additional benefit from allowing mild coordination between the trials through employing shared randomness, e.g. a shared random seed between the trials. Such a phenomenon has been observed before in the distributed testing literature \cite{acharya_inference_2020,acharya2020distributed,szabo2022optimal_IEEE,szabo2023optimal}, which forms the basis of our analysis below.

We consider the following variation on Assumption \ref{assumption:local_randomness_statistic}, where the key difference is that the source of randomness is allowed to be shared between the $m$ trials.
\begin{assumption}\label{assumption:public_randomness_statistic}
For functions $f_j: \R^d \times \R \to \R$ and a random variable $U$ which is independent of the data $X$, the $j$-th test statistic $S^{(j)} = f_j(X^{(j)},U)$ satisfies $\E_0 |S^{(j)}| \leq M$ for some $M>0$ and all $j=1,\dots,m$.
\end{assumption}
\noindent Test statistics satisfying this assumption shall be referred to as shared randomness (or public coin) protocols. 

The theorem below establishes the optimal rate when coordination through shared randomness is allowed. When the number of trials is small compared to the dimension (i.e. $m \lesssim d / \log m$), there is no difference between protocols that coordinate using shared randomness or those without coordination. In fact, the optimal rate ($\rho^2 \asymp \sqrt{d}/(\sqrt{m}n)$) in this case is reached by the test \eqref{test:large:d} or the ones below it, which do not employ shared randomness. However, when the number of trials is large compared to the dimension (i.e. $m\gtrsim d$), the testing rate substantially improves in the shared randomness protocols.

{\begin{theorem}[]\label{thm:public_coin_lb}
Let $S^{(1)},\dots,S^{(m)}$, $C_m$ and $T_{\alpha}$ satisfy Assumptions  \ref{assumption:lipschitz_type_C_m}--\ref{assumption:public_randomness_statistic}. Then there exists a constant $c> 0$ depending only on $L$, $p,q$ and $M$, such that if
\begin{equation}\label{eq:rho_condition_lb_public_coin}
\rho^2 \leq c  \frac{\left(\sqrt{m}\wedge \sqrt{\frac{d}{\log(m)}} \right) \sqrt{d}}{m n},
\end{equation}
 it holds that $\underset{f \in H_\rho}{\sup} \P_f \left( T_{\alpha} = 0 \right) > 2/3$ for all $n,m,d\in \N$ and any level $\alpha \in (0,0.1]$. 
 
At the same time, for all $\alpha,\beta\in(0,1)$ there exists a constant $C_{\alpha,\beta} > 0$ depending only on $\beta$, $L$, $p,q$, the function $\alpha \mapsto \kappa_\alpha$ and $M$, such that if 
\begin{equation}\label{eq:rho_condition_ub_private_coin}
\rho^2 \geq C_{\alpha,\beta}  \frac{\left(\sqrt{m}\wedge \sqrt{d} \right) \sqrt{d}}{m n}
\end{equation}
it holds that $\underset{f \in H_\rho}{\sup} \P_f \left( T_{\alpha} = 0 \right) \leq \beta$ for some test $T_{\alpha}$ of level $\alpha$ satisfying Assumptions  \ref{assumption:lipschitz_type_C_m}--\ref{assumption:public_randomness_statistic}. 
\end{theorem}}

\begin{remark}
Similarly to Theorem \ref{thm:private_coin_lb} the choice of ranges $0 < \alpha \leq 0.1$ and $\beta=2/3$ in the lower bound result is arbitrary, other choices are also possible as presented in the proof.
\end{remark}

A shared randomness method that attains the rate in \eqref{eq:rho_condition_ub_private_coin} is given next. Consider drawing an orthonormal $d \times d$ matrix $U$ taking values from the uniform measure on such matrices. As a test statistic, each trial computes $(\sqrt{n} U X^{(j)})_1$, which is a $N(0,1)$ random variable under the null hypothesis. A level $\alpha \in (0,1)$ meta-level test is then given by combining the local test statistics as
\begin{equation}\label{eq:pub_coin_test}
T_\alpha := \mathbbm{1}\left\{ 
\big| \frac{1}{\sqrt{m}} \underset{j=1}{\overset{m}{\sum}} (\sqrt{n} U X^{(j)})_1   \big|  \geq \Phi^{-1}(1 - \alpha/2) \right\},
\end{equation}
where $\Phi$ is the standard Gaussian CDF. {The core idea here is that for each trial, the same $1$-dimensional projection of the $d$
-dimensional data is computed, where the projection is taken uniformly at random and the test is conducted along the projected direction.} The above method corresponds to Stouffer's method for the p-values $p^{(j)} = \Phi(\sqrt{n} (U X^{(j)})_1)$ for $j=1,\dots,m$. Lemma \ref{lem:pub_coin_coordinate_strat_rate} in the appendix shows that the above test attains a small Type II error probability whenever $\rho^2 \gtrsim d/(mn)$.

\section{Examples for meta-{analysis} methods}\label{sec:examples_worked_out}

Combinations of independent test statistics that fall into the framework of Assumptions \ref{assumption:local_randomness_statistic}-- \ref{assumption:public_randomness_statistic} are subject to the rate optimality theory established by the main theorems in Section \ref{sec:main_results}. In this section, we look into common methods for combining p-values, e-values and other test-statistics, as mentioned in the introduction.

When the distribution under the null hypothesis of the test statistics are known, certain combinations are natural. For example, the sum of normal or chi-square test statistics is again normal or chi-square distributed, respectively. Similarly, voting based mechanisms typically rely on summing Bernoulli random variables. It is easy to see that these and similar combinations methods fall into the framework of Assumptions \ref{assumption:local_randomness_statistic}--\ref{assumption:public_randomness_statistic}. 

For more specific test statistics, such as p-values or e-values, many general combination methods have been introduced in the literature. We cover some of the most prominent combination approaches for p-values and e-values in Section \ref{ssec:p-values_combined} and Section \ref{ssec:e-values_combined}, respectively. The list of methods is certainly non-exhaustive and many more combination methods exist, but they serve as context for the range of techniques covered by our general theory. Our main results allow establishing lower bound rates for the ones listed below, whilst in Sections \ref{ssec:upperbound_private_coin} and \ref{ssec:public_coin} attainability of these rates by some of the listed methods was exhibited.

\subsection{Combinations of p-values}\label{ssec:p-values_combined}

If $p^{(1)},\dots,p^{(m)}$ are p-values obtained from $m$ independent test statistics concerning the same hypothesis, then under the null $p^{(j)}\sim^{iid}U(0,1)$. One can aim to combine the $m$ p-values to form a test $T_{\alpha} \equiv T_{\alpha}(p^{(1)},\dots,p^{(m)})$ with Type I error probability $\alpha$, which hopefully has higher power than a test based on one of the individual p-values. Below we list standard methods in the literature.

\begin{itemize}
\item Fisher's method \cite{fisher1992statistical}. Because the variables $- 2\log p^{(j)}$'s are iid $\chi_2^2$-distributed under the null hypothesis, their sum follows a $\chi_{2m}^2$-distribution. Therefore the combination method $ {\sum}_{j=1}^m -2  \log p^{(j)} $ results in a $\chi_{2m}^2$ distributed random variable, and the corresponding quantile function provides level-$\alpha$ one-sided tests at the meta-level.

\item Similar flavour to Fisher's method are the combinations ${\sum}_{j=1}^m - \log (1-p^{(j)})$ (Pearson's method \cite{pearson1934new}),  $\sum_{j=1}^{m} - \log p^{(j)}(1-p^{(j)})$ (the logit method / Mudholkar and George method \cite{mudholkar1977logit}) and $m^{-1/2}\sum_{j=1}^m (p^{(j)} - 1/2)$ (Edgington's method \cite{edgington1972additive}). 

\item Order-based methods such as Tippett's method \cite{tippett1941methods} based on $  \min \{ p^{(1)},\dots,p^{(m)} \} \overset{H_0}{\sim} \text{Beta}(1,m)$.

\item Methods based on inverse CDF's, such as by Stouffer et. al \cite{stouffer1949american}  based on $m^{-1/2} {\sum}_{j=1}^m  \Phi^{-1}(p^{(j)}) \sim \text{N}(0,1)$ under the null hypothesis. 

\item Generalized averages as considered in \cite{vovk2020combining},
$T_{\alpha} = \mathbbm{1} \left\{ a_{r,m} M_{r,m}(p^{(1)},\dots, p^{(m)}) \leq \alpha \right\}$, where  $M_{r,m}(p^{(1)},\dots, p^{(m)})$ equals $\big( m^{-1} {\sum}_{j=1}^m (p^{(j)})^r  \big)^{1/r}$ for $r \in \R \setminus \{ 0 \}$, the geometric mean, minimimum (i.e. Tippett's method) and maximum for $r = 0$, $r \to -\infty$, and $r \to \infty$, respectively. For $r \in \{ -\infty\}\cup [1/(m-1), \infty]$, $a_{r,m}$ can be taken to obtain precisely level $\alpha$ tests (i.e. $\P_0 T_\alpha = \alpha$). { We note that this means that canonical multiple testing methods (see e.g. \cite{familywise-error-rate}) such as Bonferroni's correction (which corresponds with taking as $M_{r,m}$ the minimum and $a_{r,m} = m$) also fall within our framework.}
\end{itemize}

Lemma \ref{lem:p_value_combs} below shows that all the methods mentioned above fall into the framework of Assumptions \ref{assumption:local_randomness_statistic}--\ref{assumption:public_randomness_statistic}. This means that the error rate lower bounds of Theorem \ref{thm:private_coin_lb} and Theorem \ref{thm:public_coin_lb} respectively, apply to the p-value combination methods listed above. That is, one cannot attain a better separation rate when considering the worst case Type II error probability for the alternative hypothesis in \eqref{eq : testing risk}, with any of the p-value combination methods listed above. Whether Assumption \ref{assumption:local_randomness_statistic} or \ref{assumption:public_randomness_statistic} applies depends on whether shared randomness is used in generating the p-values. To confirm that Assumptions \ref{assumption:test_form} and \ref{assumption:lipschitz_type_C_m} apply to tests based on the combined p-values, some algebra is needed. The proof of the lemma is deferred to the appendix.

\begin{lemma}\label{lem:p_value_combs}
Consider p-values $p^{(1)},\dots,p^{(m)}$, where each $p^{(j)}$ depends on the local data $X^{(j)}$ and possibly local randomness that is independent of the data. For each of the combination methods for p-values mentioned above and corresponding test $T_\alpha$ of level $\alpha \in (0,1)$, the conclusions of Theorem \ref{thm:private_coin_lb} holds. 
\end{lemma}

We remark that the p-values are obtained using shared randomness (i.e. in the sense of Assumption \ref{assumption:local_randomness_statistic}), the lower bound rate of Theorem \ref{thm:public_coin_lb} applies. Furthermore, as exhibited in Sections \ref{ssec:upperbound_private_coin} and \ref{ssec:public_coin}, for p-values corresponding to well chosen test statistics, these combination methods can achieve the theoretical limits established in Theorems \ref{thm:private_coin_lb} and \ref{thm:public_coin_lb}, respectively.

\subsection{Combining e-values}\label{ssec:e-values_combined}

An \emph{e-value} is a nonnegative random variable $E$ such that $\sup_{\P_0\in H_0}\P_0 E \leq 1$. The \emph{threshold test corresponding to $E$ of level $\alpha$} is $\mathbbm{1} \{ E \geq \alpha^{-1} \}$. This test yields a so called strict p-value; for $\P_0 \in H_0$ we have $\P_0( E \geq \alpha^{-1}) \leq \alpha$ by Markov's inequality.

E-values lend themselves for combining outcomes of independent studies for two main reasons. First, they are easy to combine, see Section 4 in \cite{vovk2021values} for an indepth  discussion of specific combination functions for independent e-values. Second, they are robust to misspecification and offer optional stopping/continuation guarantees \cite{IEEE_Safe_Testing}. Common examples of e-values are Bayes factors and likelihood ratios{, which are nonnegative and have expectation equal to $1$ in the case of a simple null hypothesis such as considered in this article.}

Several combination methods (e-merging functions) were proposed in the literature. For instance, the product of independent e-values is also again an e-value. This was shown to weakly dominate any other combination of independent e-values in the sense that $\Pi_{j=1}^{m} E^{(j)} \geq C_m(E)$, for any  $E = (E^{(j)}) \in [1,\infty)^m$ and $E \mapsto C_m(E)$ such that $C_m(E^{(1)},\dots,E^{(m)})$ is an e-value for any independent e-values $E^{(1)},\dots,E^{(m)}$, see \cite{vovk2021values}. Similarly, the average of e-values is again an e-value. The product and the average are \emph{admissible} in the sense that there is no e-merging function that strictly dominates them on $[0,\infty]^m$. The lemma below shows that these two, arguably most prominent e-value combination methods fulfill Assumptions \ref{assumption:local_randomness_statistic}-- \ref{assumption:public_randomness_statistic} and hence the lower bounds derived in Theorems \ref{thm:private_coin_lb} and   \ref{thm:public_coin_lb} apply.

\begin{lemma}\label{lem:e_value_combs}
Consider e-values $E^{(1)},\dots,E^{(m)}$, where each $E^{(j)}$ depends on the local data $X^{(j)}$ and possibly local randomness that is independent of the data. Let $C_m:\R^m \to \R$ correspond to either the average or the product and let $T_\alpha$ be the corresponding threshold test of level $\alpha \in (0,1)$,
\begin{equation*}
T_\alpha = \mathbbm{1} \left\{ C_m(E^{(1)},\dots,E^{(m)}) \geq \alpha^{-1} \right\}.
\end{equation*}
If $C_m$ is the product, assume in addition that $\E_0 | \log E^{(j)} | $ is uniformly bounded. Then, the conclusion of Theorem \ref{thm:private_coin_lb} holds. In case the e-values are generated using shared randomness, then Theorem \ref{thm:public_coin_lb} applies.
\end{lemma}

\section{Simulations}\label{sec:sim}

In this section, we investigate the numerical performance of the testing strategies outlined in Section \ref{ssec:upperbound_private_coin} on synthetic data sets. We compare the tests based on their receiver operating characteristic (ROC) curve. For a range of significance levels we compute for each tests  the ``true positive rate'' (TPR) and ``false positive rates'' or (FPR), i.e.
the fraction of the simulation runs in which the test correctly identifies the underlying signal, falsely rejects the null hypothesis, respectively. Plotting the TPR against the FPR (both given as a function of the significance level) provides us the ROC curve, visualizing the  diagnostic ability of the test.

\begin{figure}[ht]\label{fig:roc_curves}
\includegraphics[width=0.9\textwidth]{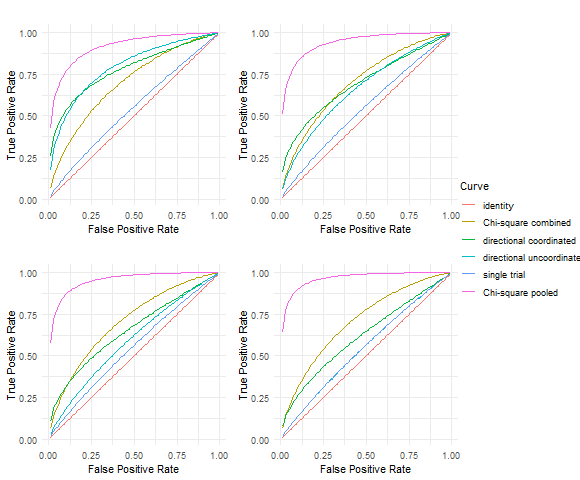}
\caption{ROC curves for different values of $d$, whilst keeping $m=20$, $n=30$, $\rho^2 =  \sqrt{d}/(4n)$. From left to right, top to bottom: $d=2$, $d=5$, $d=10$, $d=20$.}
\end{figure}

In our simulations we set $m=20$, $n=30$, let $d$ range from $2$ to $20$ and take $\rho^2 =  \sqrt{d}/(4n)$. This value of $\rho^2$ corresponds to a signal that is almost indistinguishable {from} noise using just a single trial, whilst consistently detectable if the data were to be pooled with $m \approx 20$ (which increases the signal size to noise ratio effectively by a factor $\sqrt{20} > 4$. For each level $\alpha \in \{ 0.01, 0.02,\dots, 0.99 \}$ we compute the power for different combination strategies $100$ times, each time drawing a different $f \in \R^d$ with $\| f \|_2 = \rho$ according to $f_i = d^{-1/2} \rho R_i$ and $R_i$ iid Rademacher random variables for $i=1,\dots,d$. As combination strategies, we compare the strategies \eqref{test:large:d}, \eqref{eq:pub_coin_test} and \eqref{eq:private_coin_test} from Section \ref{sec:main_results}, which are called ``chi-square combined'', ``coordinated directional'' and ``uncoordinated directional'' in the legend of Figure \ref{fig:roc_curves}. In addition, we display the ROC curves for the chi-square test based on pooled data (``chi-square pooled'') and that of a single trial (``single trial'').

We make the following observations, in line with our theoretical findings. The {meta-analysis} methods based on combining the locally optimal chi-squared test statistics (yellow curves) substantially out-performed the chi-squared test statistics based on a single trial (blue curve), but was substantially worse than the chi-square test based on the pooled data (pink curve). Second note that the large dimensional case ($d=10$ and $d=20$) the best strategy is indeed to combine the local chi-square statistics (yellow curve), while in the low dimensional setting ($d=2$) it is more advantageous to combine the directional test statistics $X_i^{(j)}$ (blue curve). Finally, note that allowing coordination between the trials by a shared randomness protocol can result in improved performance (green curve) compared to the independent experiments (blue curve). In fact this approach provides the best {meta-analysis} method in the small dimensional setting (e.g. $d=2$ and $d=5$ for small $\alpha$, which is the most interesting case). 

{In the appendix, Section \ref{sec:additional_sims}, we explore eight additional simulation settings, where we consider larger values of $d$ and $m$. Whilst these simulations do not reveal additional phenomena to the ones observed in Figure \ref{fig:roc_curves}, they do give insight into the relative performance of the testing methods for different values of $d$ and $m$.}

\section{Discussion}\label{sec:discussion}

We briefly summarize our main contributions and discuss possible extensions and research directions. First, by establishing a connection between {meta-analysis} and distributed learning under communication constraints, we have provided a unified, theoretical framework for evaluating the behaviour of standard meta-analysis techniques. {In our analysis, we considered the many normal means model, but these results can be extended to other more complex models as well, building on the connection with distributed computation. For example, minimax estimation rates under communication constraints were derived for other parametric models \cite{zhang2013information},  density estimation \cite{barnes2020lower}, signal-in-Gaussian-white-noise \cite{zhu2018distributed,szabo2022distributed,cai2022distributed}, nonparametric regression \cite{szabo_adaptive_2019} and in abstract settings \cite{zaman2022distributed} including binary and Poisson regression, drift estimation, and more. The normal means model allows for a tractable analysis, but results in this model are known to extend to more complicated models, such as discrete density testing (see e.g. \cite{carter2002multinomial}). With the due technical work, our results are expected to translate to these settings as well, but we leave this for future endeavor.}

In the normal means model we show that by combining the locally optimal chi-square statistics at a meta-level one can gain a factor of $\sqrt{m}$ compared to using a single trial. Nevertheless, regardless of the choice of the combination method, a factor of $\sqrt{m}$ is lost compared to the scenario when all data from all trials are at our disposal. This loss is clearly visible even in small sample sizes, dimensions and trial numbers, as demonstrated in our numerical analysis, as can be seen in the corresponding ROC curves. For more complex models, such a numerical study can be a first step to quantify the efficiency of the {meta-analysis} method. We have also shown that in the small dimension - large number of trials setting combining the locally optimal chi-square statistics (or any rotationally invariant statistics for that matter) results in information loss and sub-optimal accuracy. In this case, better rates can be attained by test statistics based on the direction of the observations combined at the meta-level. {In practice, one often cannot choose which test statistics can be obtained from independent trials. In such cases, the $\sqrt{m}$-factor loss in the case of e.g. rotationally invariant test statistics is of interest when considering power calculations. Meta-analysis approaches based on directional test statistics are designed for scenarios where individual datasets are not centrally collected, but there is some level of coordination among experimenters.}

{The assumption throughout the paper of homogeneity between the trials (i.e. each trial consisting of the same number of observations) simplifies the presentation, but the results can be extended to cases where the number of observations in each trial differ by constant factors. Situations where the number of observations differs greatly (e.g. $k \ll m$ trials have as much observations as the other $m-k$ trials combined) are certainly of interest, but beyond the scope of this paper.}

\pagebreak 



\textbf{{Acknowledgements}:} Co-funded by the European Union (ERC, BigBayesUQ, project number: 101041064). Views and opinions expressed are however those of the author(s) only and do not necessarily reflect those of the European Union or the European Research Council. Neither the European Union nor the granting authority can be held responsible for them. This research was partially funded by a Spinoza
grant of the Dutch Research Council (NWO).

\bibliographystyle{acm}
\bibliography{references}

\pagebreak 

\setcounter{section}{0}
\setcounter{equation}{0}
\renewcommand{\theequation}{S.\arabic{equation}}
\renewcommand{\thesection}{\Alph{section}}

\title{Supplementary material to ``optimal testing using combined test statistics across independent studies''}

\begin{abstract}

In this supplement, we present the detailed proofs for the main theorems in the paper ``Optimal testing using combined test statistics across independent studies''.

\end{abstract}

\section{Appendix}

The proofs of the main theorems (Theorem \ref{thm:private_coin_lb}, \ref{thm:private_coin_ub} and \ref{thm:public_coin_lb}) are divided over the subsections as follows. In Section \ref{ssec:main_theorem_proof}, the lower bounds of Theorem \ref{thm:private_coin_lb} and \ref{thm:public_coin_lb} are proven. Auxiliary lemmas for the proof of the lower bounds are proven in \ref{ssec:auxilliary_lemmas}. The attainability of the lower bound rates are given in Lemmas \ref{lem:sqrtm_rate_attained}, \ref{lem:priv_coin_coordinate_strat_rate} and \ref{lem:pub_coin_coordinate_strat_rate} in Section \ref{ssec:rates_attainable}. In Section \ref{ssec:proofs_p_e_values}, Lemmas \ref{lem:p_value_combs} and \ref{lem:e_value_combs} are proven.

\subsection{Proof of the lower bounds (Theorems \ref{thm:private_coin_lb} and \ref{thm:public_coin_lb})}\label{ssec:main_theorem_proof}

The proof is based around the following idea. If $C_m$ satisfies the continuity condition of Assumption \ref{assumption:lipschitz_type_C_m}, it implies $C_m(S^{(1)},\dots,S^{(m)})$ should not change to much if the statistics $S^{(1)},\dots,S^{(m)}$ are replaced by finite bit approximations. If $b$ is the number of bits used for the approximation of $S^{(j)}$, we should be able to get an approximation with accuracy of the order $2^{-b}$ through e.g. binary expansion. Since $C_m$ and consequently the test based on $C_m$ do not change (much) from passing to a finite bit approximation, tools and results from testing under bit-constrained communication apply, which finally yield the theorems.

\begin{proof}
We prove the statement for any $\alpha \in (0,1/10]$. Since $\alpha \mapsto \kappa_\alpha$ is strictly decreasing, $\kappa_{1/8} < \kappa_{1/10} \leq \kappa_{\alpha}$ holds  for any $\alpha \in (0,1/10]$. Take $0<\epsilon<\frac{1}{2}(\kappa_{1/10}-\kappa_{1/8})$. Then $|x - \kappa_\alpha| \leq \epsilon$ implies $x \geq \kappa_{1/8}$, which by the definition of the quantile function provides
\begin{equation}\label{eq:C_m_sufficiently_anticoncentrated}
\P_0 \left( |C_m(S) - \kappa_\alpha| \leq 2\epsilon \right) \leq 1/8.
\end{equation}
By Lemma \ref{lem:binary_expansion_lemma}, there exist $B^{(j)}$-bit binary approximations $\tilde S^{(j)}$ such that
\begin{equation}\label{eq:tilde_S_approx_well}
| S^{(j)} - \tilde S^{(j)}| \leq \left(\frac{\epsilon^{1/q}}{L^{1/q}m}\right)^{1/p}
\end{equation}
and 
\begin{equation}\label{eq:bits_expectation_bound}
\E_0 B^{(j)} \leq \E_0 \log_2 (|S^{(j)}|)\vee 0 - \frac{1}{p} \log \left( \frac{\epsilon^{1/q}}{L^{1/q}m} \right) + 3.
\end{equation}
Write $\tilde{S} = (\tilde{S}^{(1)}, \dots, \tilde{S}^{(m)})$. By combining Assumption \ref{assumption:lipschitz_type_C_m} with \eqref{eq:tilde_S_approx_well},
\begin{equation*}
| C_m(S) - C_m(\tilde S) | \leq \epsilon.
\end{equation*}
Consequently,
\begin{align*}
\cR(T_\alpha, H_\rho) &\geq \P_0 \left( C_m(\tilde S) - | C_m(S) - C_m(\tilde S) | \geq \kappa_\alpha \right) \\ & \;\;\;\;+ \underset{f \in H_\rho}{\sup} \P_f \left( C_m(\tilde S)  \leq \kappa_\alpha - |C_m(S) - C_m(\tilde S)| \right)
\\ &\geq \P_0 \left( C_m(\tilde S) \geq \kappa_\alpha + \epsilon \right) + \underset{f \in H_\rho}{\sup} \P_f \left( C_m(\tilde S)  \leq \kappa_\alpha - \epsilon \right).
\end{align*}
Define the test
\begin{equation*}
T'_\alpha :=\mathbbm{1} \left \{ C_m(\tilde S) > \kappa_\alpha - \epsilon \right\}.
\end{equation*}
Since
\begin{equation*}
\P_0 \left( C_m(\tilde S) \geq \kappa_\alpha + \epsilon \right) = \P_0 \left( C_m(\tilde S) > \kappa_\alpha - \epsilon \right) - \P_0 \left( - \epsilon \leq C_m(\tilde S) - \kappa_\alpha \leq   \epsilon \right),
\end{equation*}
the second last display can now be written as
\begin{equation*}
\cR(T',H_\rho) - \P_0 \left( |C_m(\tilde S) - \kappa_\alpha| \leq \epsilon \right).
\end{equation*}
Applying \eqref{eq:C_m_approximation_assumption} again, using the reverse triangle inequality and \eqref{eq:C_m_sufficiently_anticoncentrated}, we obtain
\begin{equation*}
\P_0 \left( |C_m(\tilde S) - \kappa_\alpha| \leq \epsilon \right) \leq \P_0 \left( |C_m(S) - \kappa_\alpha| \leq 2\epsilon \right) \leq 1/8.
\end{equation*}
It suffices to show that for $\rho$ satisfying \eqref{eq:rho_condition_private_coin_lb} in the case of Theorem \ref{thm:private_coin_lb} or $\rho$ satisfying \eqref{eq:rho_condition_lb_public_coin} in case of Theorem \ref{thm:public_coin_lb} for a small enough $c>0$, we have 
\begin{equation}\label{eq:to_show_risk_lower_bound}
\cR(T',H_\rho) \geq 7/8.
\end{equation}
This follows from Lemma \ref{lem:workhorse}, where it is left to verify that 
\begin{equation}\label{eq:to_show_bit_bound}
\underset{j=1}{\overset{m}{\sum}} {d \wedge \E_0 {B}^{(j)}} \lesssim m (d \wedge (1 \vee \log m))
\end{equation}
for a constant independent of $d,n,m$ and $c>0$. By \eqref{eq:bits_expectation_bound} and $\E_0 |S^{(j)}|\leq M$ for some constant $M>0$ for $j=1,\dots,m$ (following from Assumption \ref{assumption:local_randomness_statistic} or \ref{assumption:public_randomness_statistic}), we obtain that $\underset{j=1}{\overset{m}{\sum}} {d \wedge\E_0 {B}^{(j)}}$ is bounded by 
\begin{equation*}
  m \left( d \bigwedge \left( 
 \log_2 (1+M)   + 3 - \frac{1}{p} \log \left( \frac{\epsilon^{1/q}}{L^{1/q}m} \right)  \right)\right),
\end{equation*}
from which \eqref{eq:to_show_bit_bound} follows. Putting things together, we now have that for $c > 0$ small enough we obtain \eqref{eq:to_show_risk_lower_bound}, from which we conclude that \eqref{eq:to_show_risk_lower_bound} holds and the proof of the theorems is concluded.
\end{proof}

\subsection{Auxiliary lemmas to the lower bound theorems}\label{ssec:auxilliary_lemmas}

As a first tool, we introduce finite bit approximations of real numbers through their binary expansion. Consider the binary expansion of $x \in \R$; i.e. there exist digits $a_k(x),\dots,a_1(x),a_0(x) \in \{0,1\}$ for a $k_x \equiv k \in \N \cup \{ 0 \}$ and $(b_i(x))_{i \in \N} \in \{0,1\}^\N$ such that 
\begin{equation}\label{eq:B_bit_binary_expansion}
x = \text{sign}(x) \left(\underset{i=0}{\overset{k}{\sum}} 2^i a_i(x) + \underset{i=1}{\overset{\infty}{\sum}} 2^{-i} b_i(x)  \right)
\end{equation}
with $k$ the largest element in $\N \cup \{ 0 \}$ such that $ 2^k - 1 \leq |x|$. We now define $\tilde{x}_B$ to be the $B$-bit binary expansion giving the smallest approximation error in absolute value, where the first bit encodes $\text{sign}(x)$. That is, for $B \geq k + 2$, we have
\begin{equation}\label{eq:binary_approximation_accuracy_general}
| x - \tilde{x}_B | \leq \underset{i=B - k -1}{\overset{\infty}{\sum}} 2^{-i} b_i(x).
\end{equation}
The following is well known, we exhibit its proof for completeness.

\begin{lemma}\label{lem:binary_expansion_lemma}
Let $V$ be a random variable with a first moment. Given $1>\epsilon > 0$, let $B_\epsilon \equiv B$ denote the number of bits required such that
\begin{equation}\label{eq:binary_accuracy}
|V - \tilde{V}_{B_\epsilon}| \leq \epsilon.
\end{equation}
It holds that
\begin{equation*}
\E B \leq \E {\log_2(|V|)\vee 0 + 1} + \log_2(1/\epsilon)  + 2.
\end{equation*}
\end{lemma}
\begin{proof}
If $|V| < 1$, we have that
\begin{equation*}
|V - \tilde{V}_B| \leq \underset{i=B -1}{\overset{\infty}{\sum}} 2^{-i} b_i(V).
\end{equation*}
So in the case that $|V| \leq 1$, since $ b_i(V) \in \{0, 1\}$, for \eqref{eq:binary_accuracy} to hold it suffices that $B \geq \log_2(1/\epsilon) + 2$. Let $B'$ denote the amount of bits required to obtain $| V - \tilde{V}_{B'} | \leq 1$. When $2^k \leq |V| < 2^{k+1}$, it holds that $B' \leq k+1$. Using Markov's inequality,
\begin{align*}
\E B' &= \E B' \underset{k=0}{\overset{\infty}{\sum}} \mathbbm{1}\left\{ 2^k \leq |V| < 2^{k+1} \right\} \\
&\leq \E  \underset{k=0}{\overset{\infty}{\sum}}  ({k+1}) \mathbbm{1}\left\{ k \leq \log_2(|V|) < {k+1} \right\} \leq  \E {\log_2(|V|)\vee 0 + 1}.
\end{align*}
In conclusion, $\E B \leq \E {\log_2(|V|)\vee 0 + 1} + \log_2(1/\epsilon)  + 2$.
\end{proof}

For the lemmas below, we introduce the following notation. Let $\pi$ be a probability distribution on $\R^d$. Write $\P_\pi := \int \P_f d\pi(f)$ for the mixture distribution, where $\P_f$ denotes the joint distribution on $X$, $U$ and $S$. Let $F$ denote the draw from $\pi$. Let $\P_f^{\tilde{S}}$ denote the forward measure induced on the random variable $\tilde{S}$ and let $L_\pi^{\tilde{S}}$ denote the likelihood ratio of the mixture distribution and $\P_0$, ie
\begin{equation}\label{eq:likelihood_S}
L_\pi^{\tilde{S}} = \int \frac{d\P^{\tilde{S}}_f}{d\P^{\tilde{S}}_0} d\pi(f).
\end{equation}
 Because of the Markov chain structure of $F \to (X,U) \to S$ and the independence between $U$ and $X$, the joint distribution of $(X,U,S)$ under the mixture disintegrates as
\begin{equation}\label{eq:measure_factorizes}
d\P^{X,U,S}_\pi(x,u,s) = \int d\P^{S|(X,U)}(s) d\P_f^{X}(x) d\P^{U}(u) d\pi(f)
\end{equation}
where $\P^{U}$ is the marginal distribution of $U$. For the likelihood ratio conditionally on $U=u$, we shall write
\begin{equation}\label{eq:likelihood_S_given_U}
L_\pi^{\tilde{S}|U=u} = \int \frac{d\P^{\tilde{S}|U=u}_f}{d\P^{\tilde{S}|U=u}_0} d\pi(f).
\end{equation}
Furthermore, by the independence of the statistics given $U$, 
\begin{equation}\label{eq:measure_factorizes_locally}
d\P^{S|(X,U)} = \underset{j=1}{\overset{m}{\bigotimes}}   d\P^{S^{(j)}|(X^{(j)},U)}.
\end{equation}
Let $\tilde S^{(j)}$ denote the $B^{(j)}$-bit binary approximations to $S^{(j)}$ such that \eqref{eq:tilde_S_approx_well} holds. Note that the above displays are true for the random variable $\tilde{S} = (\tilde S^{(1)},\dots,\tilde S^{(m)})$ in place of $S$ since $F \to (X,U) \to S \to \tilde{S}$ forms a Markov chain as well. The following lemma allows us to bound the chi-square divergence between the forward measure for $\tilde{S}$, which we will denote by $\P_\pi^{\tilde{S}}$ and $\P_0^{\tilde{S}}$.

The following lemma lower bounds the worst-case risk for any test $T'$ depending only on $\tilde{S}$, the binary approximation of $S$ as in \eqref{eq:tilde_S_approx_well}.
\begin{lemma}\label{lem:workhorse}
Let $T'$ be a test depending only on $\tilde{S}$ taking values in $\R^m$, satisfying \eqref{eq:measure_factorizes} and where $\tilde{S}^{(j)}$ allows for an exact ${B}^{(j)}$-bit binary expansion as in \eqref{eq:B_bit_binary_expansion}, with $\E_0 {B}^{(j)} < \infty$ for $j=1,\dots,m$.

There exists $c > 0$ independent of $n,m$ and $d$ such that
\begin{equation*}
\cR(T',H_\rho) \geq 7/8
\end{equation*}
for all $n,m,d \in \N$ whenever 
\begin{equation}\label{eq:trace_bit_bound}
  \underset{j=1}{\overset{m}{\sum}} {d \wedge \E_0 {B}^{(j)}} \lesssim m (d \wedge \log m)
\end{equation}
in addition to
\begin{equation}\label{eq:appendix_rho_bound_public_coin}
\rho^2 \leq c  \frac{(\sqrt{m}\wedge \frac{d}{\log(m)}) \sqrt{d}}{m n},
\end{equation}
if $\tilde{S}$ is generated using public randomness, or
\begin{equation}\label{eq:appendix_rho_bound_private_coin}
\rho^2 \leq c \frac{(\sqrt{m}\wedge \sqrt{\frac{d}{\log(m)}}) \sqrt{d}}{m n},
\end{equation}
in case $\tilde{S}$ is generated using only local randomness.
\end{lemma}
\begin{proof}
Consider a probability distribution $\pi$ on $\R^d$ and $L_\pi^{\tilde{S}}$ as in \eqref{eq:likelihood_S}. 
Consider the set
\begin{equation*}
D := \left\{ u : \underset{j=1}{\overset{m}{\sum}}  d \wedge \E_0[ B^{(j)} | U=u] \leq 64 \underset{j=1}{\overset{m}{\sum}}  d \wedge  \E_0 B^{(j)} \right\},
\end{equation*}
whose complement, {$D^c$}, has $\P^U$-mass less than or equal to $1/64$ by Markov's inequality and $\E^U (d \wedge \E_0[ B^{(j)} | U]) \leq d \wedge \E_0 B^{(j)}$. By conditioning on $U$ (writing $\P^{|U=u}_0 := \P_0(\cdot | U=u)$),
\begin{align*}
\cR(T',H_\rho) &\geq \P_0 T' + \P_\pi (1 - T') - \pi(f \notin H_\rho) \\
&\geq \int \left(\P_0^{|U=u} (T') + \P_\pi^{|U=u} (1 - T') \right)\mathbbm{1}_D(u) d\P^{U}(u)- \pi(f \notin H_\rho).
\end{align*}
Since $0\leq T' \leq 1$ and $L_\pi^{\tilde{S}} \geq 0$, for all $0 < \gamma < 1$,
\begin{align*}
            \P_0^{|U=u} (T') + \P_\pi^{|U=u} (1 - T')    &\geq \P_0^{|U=u} \left(\gamma T' + L_\pi^{\tilde{S}|U=u} (1 - T') \mathbbm{1}\left\{  L_\pi^{\tilde{S}|U=u} > \gamma \right\}\right)  \\
                &\geq \gamma \P_0^{|U=u} \left(  L_\pi^{\tilde{S}|U=u} > \gamma \right) \\
                &\geq \gamma\big(1 - \P_0^{|U=u} (  |L_\pi^{\tilde{S}|U=u} - 1| \geq 1 - \gamma )\big).
\end{align*}
The probability on the right hand side of the above display can be bounded by applying Chebyshev's inequality and bounding the resulting chi-square divergence using the tools of \cite{szabo2023optimal}, in particular using Lemma 10.1 from the aforementioned paper. This lemma applies if $\tilde{S}$ takes values in a space of finite, fixed cardinality. 

Define $B^* = \underset{j=1}{\overset{m}{\sum}} {64 \E_0 |B^{(j)}|}$ and the event
\begin{equation*}
A := \left\{ \underset{j=1}{\overset{m}{\sum}} B^{(j)} \leq B^*  \right\},
\end{equation*}
so that $A^c$ by Markov's inequality occurs with $\P_0$-probability less than $1/64$. 

Let $\breve{S}^{(j)}$ be the $\breve{B}^{(j)}:= {B}^{(j)} \wedge B^*$ binary approximation of $\tilde{S}^{(j)}$ and note that on the event $A$, $\breve{S}^{(j)} = \tilde{S}^{(j)}$. We have
\begin{align*}
\int \P_0^{|U=u} \left(  |L_\pi^{\tilde{S}|U=u} - 1| \geq 1 - \gamma \right) \mathbbm{1}_D(u) d\P^{U}(u) &\leq  \\ \int \P_0^{|U=u} \left( \left\{ |L_\pi^{\tilde{S}|U=u} - 1| \geq 1 - \gamma \right\} \cap A \right) \mathbbm{1}_D(u) d\P^{U}(u) + \P_0(A^c) &\leq \\
 \int \P_0^{|U=u} \left( |L_\pi^{\breve{S}|U=u} - 1| \geq 1 - \gamma  \right) \mathbbm{1}_D(u) d\P^{U}(u) + 1/64, \, &
\end{align*}
where $\breve{S} = (\breve{S}^{(1)},\dots,\breve{S}^{(m)})$. Using \eqref{eq:measure_factorizes} and Chebyshev's inequality, it suffices to show that on the event $D$, $\E_0^{|U=u} |L_\pi^{\breve{S}|U=u} - 1|^2$ is smaller than $\frac{1}{32}(1-\gamma)^2$ for $c$ small enough when $\rho$ satisfies \eqref{eq:appendix_rho_bound_public_coin} or \eqref{eq:appendix_rho_bound_private_coin}, some $\gamma \geq 5/6$ for a specific choice of $\pi$. By Lemma \ref{lem:chi-square-bound}, such a distribution $\pi$ exists, satisfying $\pi(f \notin H_\rho) \leq 1/32$, as long as $\text{Tr}(\Xi_u)$ can be sufficiently bounded, which can be done in terms of \eqref{eq:bits_expectation_bound}, as we will show next.

Let $\cS^{(j)}(b,u)$ be the space in which $\breve{S}^{(j)}|[B^{(j)} = b,U=u]$ takes values. Write 
\begin{equation*}
V_{s,u} = \E_0\left[ {X}^{(j)} \bigg| \breve S^{(j)}=s, U=u\right].
\end{equation*}
We have
\begin{align*}
\Xi^j_u &= \sum_s V_{s,u} V_{s,u}^\top \P^{}_0(\breve{S}^{(j)}=s|U=u) \\
 &= \underset{b \in \N}{\overset{}{\sum}} \P_0( \breve{B}^{(j)} = b|U=u)  \underset{s = \cS^j(b,u) }{\overset{}{\sum}} \P_0( \breve S^{(j)}=s | \breve{B}^{(j)}=b, U=u) V_{s,u} V_{s,u}^\top.
\end{align*}
By Lemma A.3 in \cite{szabo2023optimal}, the trace of the matrix
\begin{equation*}
\underset{s \in \cS^j(b,u) }{\overset{}{\sum}} \P_0 \left( \breve S^{(j)}=s | \breve{B}^{(j)}=b, U=u\right) V_{s,u} V_{s,u}^\top
\end{equation*}
is bounded by $\left(2 \log(2) \frac{b}{d} \bigwedge 1 \right) \frac{d}{n}$. By linearity of the trace operation,
\begin{align*}
\text{Trace}(\Xi^j_u) &= \underset{b \in \N}{\overset{}{\sum}} \P_0 \left( \breve{B}^{(j)} = b|U=u\right) \left(2 \log(2) \frac{b}{d} \bigwedge 1 \right) \frac{d}{n} \\ &\leq 2 \log(2) \frac{d \wedge \E_0 [\breve{B}^{(j)} | U=u]}{n}
\end{align*}
and consequently, since $\breve{B}^{(j)} \leq {B}^{(j)}$ and $u \in D$,
\begin{align*}
\text{Trace}(\sum_{j=1}^{m} \Xi^j_u) &\leq 2\log(2) n^{-1}\sum_{j=1}^{m} d \wedge \E_0 \left[  \breve{B}^{(j)} | U=u \right] \\
&\leq 128 \log(2) n^{-1} \sum_{j=1}^{m} d\wedge  \E_0 \left[  {B}^{(j)} \right].
\end{align*}
The result follows after using that $\rho^2$ satisfies \eqref{eq:appendix_rho_bound_private_coin} and \eqref{eq:appendix_rho_bound_public_coin} in the case of local or shared randomness protocols, respectively.
\end{proof}

\begin{lemma}\label{lem:chi-square-bound}
Let $L_\pi^{\breve{S}}$ be as defined through \eqref{eq:likelihood_S}, with $\breve{S}=(\breve{S}^{(1)},\dots,\breve{S}^{(m)})$ taking values in a space of finite cardinality. Let $\Xi_u = \sum_{j=1}^{m} \Xi^j_u$ with
\begin{equation}\label{eq:xi_definition}
\Xi^j_u := \E_0^{|U=u} \E_0\left[ {X}^{(j)} \bigg| \breve{S}^{(j)}, U=u\right] \E_0\left[ {X}^{(j)} \bigg| \breve{S}^{(j)}, U=u \right]^\top.
\end{equation}
Let $\rho^2$ satisfy \eqref{eq:appendix_rho_bound_public_coin} or \eqref{eq:appendix_rho_bound_private_coin}. For $c > 0$ small enough (in \eqref{eq:appendix_rho_bound_public_coin} or \eqref{eq:appendix_rho_bound_private_coin}) there exists a probability distribution $\pi$ on $\R^d$ such that
\begin{equation}\label{eq:to_show_little_mass_outside}
\pi( f \notin H_\rho ) \leq 1/32
\end{equation}
and
\begin{align}\label{eq:public_coin_chisq_divergence}
\E_0^{|U=u} |L_\pi^{\breve{S}|U=u} - 1|^2 \leq  \exp \left( C(\frac{m n^2 \rho^4}{c d} + \frac{m n^3 \rho^4}{ d^2 c}  \text{Tr}\left( \Xi_u  \right) \right) -1,
\end{align}
for a constant $C > 0$ that does not depend on $d,n,m$ or $c$. Furthermore, in case of private coin randomness ($U$ is degenerate), there exists a probability distribution $\pi$ on $\R^d$ such that \eqref{eq:to_show_little_mass_outside} is satisfied and (the sharper bound)
\begin{align}\label{eq:private_coin_chisq_divergence}
\E_0 |L_\pi^{\breve{S}} - 1|^2 \leq   \exp \left( C(\frac{m n^2 \rho^4}{c d} + \frac{n^4 \rho^4}{ d^3 c} \text{Tr}\left( \Xi_u  \right)^2 \right) - 1
\end{align}
holds for $c > 0$ small enough.
\end{lemma}

\begin{proof}
This follows from the proof of Theorem 3.1 in \cite{szabo2023optimal} (where it is important to note that in the notation of \cite{szabo2023optimal}, ``$n$'' equals ``$nm$'' in this article). For completeness, we highlight the main steps here. We start by noting that
\begin{align*}
\E_0^{|U=u} |L_\pi^{\breve{S}|U=u} - 1|^2 &= D_{\chi^2}(\P^{\breve{S}|U=u}_{0}  ; \P^{\breve{S}|U=u}_{\pi}).
\end{align*}
Let $\pi$ be a $N(0,\Gamma)$-distribution with $\Gamma \in \R^{d \times d}$. In view of the Markov chain structure (ie \eqref{eq:measure_factorizes} and \eqref{eq:measure_factorizes_locally}), the Gaussianity of $\pi$ and the fact that $\breve{B}^{(j)}$ is bounded for $j=1,\dots,m$, we obtain through following the steps corresponding to displays (34) up until (42) in Section 9 of \cite{szabo2023optimal} that the above display is bounded by
\begin{equation}\label{eq : continue to bound this display}
\underset{j=1}{\overset{m}{\Pi}} \E_0^{X^{j}|U=u} \left[ \mathscr{L}_\pi \left({X}^{j}\right)^2 \right] \cdot \int e^{\frac{n^2}{m^2} f^\top \sum_{j=1}^{m} \Xi^j_u g} d(\pi \times  \pi)(f,g)-1,
\end{equation}
where $\mathscr{L}_\pi \left({X}^{j}\right) = \int \frac{d\P^{X^{(j)}}_f}{d\P^{X^{(j)}_0}} d\pi(f)$ and we note that Lemma 10.1 applies by the boundedness of $\breve{B}^{(j)}$ and Gaussianity of $\pi$. Taking $\Gamma \in \R^{d \times d}$ equal to
\begin{equation*}
\frac{\rho}{c^{1/4} \sqrt{d}}\bar{\Gamma}
\end{equation*}
for a symmetric idempotent $d \times d$ matrix $\bar{\Gamma}$ with rank (proportional to) $d$, we obtain \eqref{eq:to_show_little_mass_outside} for $c > 0$ small enough (Lemma A.13 in \cite{szabo2023optimal}). Following the second step of the proof of Section 9 in \cite{szabo2023optimal}, in particular the steps corresponding to displays (43) and (44), we obtain that \eqref{eq : continue to bound this display} is bounded by
\begin{equation*}
\exp \left( C(\frac{m n^2 \rho^4}{c d} + \frac{n^4 \rho^4}{ d^2 c}  \text{Tr}\left( (\sqrt{\bar{\Gamma}} \Xi_u \sqrt{\bar{\Gamma}} )^2 \right) \right) - 1
\end{equation*}
for some fixed constant $C > 0$ independent of $d,n,m$ and $\rho$. The shared randomness bound of \eqref{eq:public_coin_chisq_divergence} now follows by choosing of $\bar{\Gamma} = I_d$ and using that $\text{Tr}(A^\top A) \leq \|A\| \text{Tr}(A)$ where $\| A \|$ is the operator norm of $A$, as well as by the fact that $ \Xi_u \leq \frac{m}{n}I_d$ (see Lemma A.2 in \cite{szabo2023optimal}). In case of private randomness, we can assume that $U$ is degenerate, so $\Xi_u = \Xi$ for $\P^U$-almost every $u$. The matrix $\Xi$ is positive definite and symmetric, therefore it possesses a spectral decomposition $V^\top \text{Diag}(\xi_1,\dots,\xi_d) V$. Assuming that $\xi_1 \geq \xi_2 \geq \ldots \geq \xi_d$ with corresponding eigenvectors $V = \left(\begin{matrix} v_1 & \dots & v_d \end{matrix}\right)$, let $\check{V}$ denote the $d \times \lceil d/2 \rceil$ matrix $\left(\begin{matrix} v_{\lfloor d/2 \rfloor+1} & \dots & v_d \end{matrix}\right)$. The bound of \eqref{eq:private_coin_chisq_divergence} now follows by setting $\bar{\Gamma} = \check{V} \check{V}^\top$, for a detailed computation, see page 23 of \cite{szabo2023optimal}.
\end{proof}

\subsection{Theorem concerning necessity of signs}

The theorem below tells us that in order to attain the rate of $\frac{d}{nm}$, the statistics $S^{(j)}$ need to contain at least \emph{some} information on the signs of $X^{(j)}$, in the sense that $\sqrt{d}/(\sqrt{m}n)$ is the rate that can be attained at best when $S^{(j)}$ is measurable with respect to the absolute values of the coordinates of $X^{(j)}$. This is in particular the case for statistics based on e.g. the norm $\| X^{(j)} \|_2$ or rotation invariant statistics such as the worst-case growth rate optimal e-values (see e.g. \cite{IEEE_Safe_Testing}), which consequently attain the rate $\frac{\sqrt{d}}{\sqrt{m}n}$ at best and are thus suboptimal when $d$ is small compared to $m$.

\begin{theorem}\label{thm:sign_info_necessary_for_rate}
Suppose that $S^{(j)} = f_j(X^{(j)},U)$ is such that $S^{(j)}$ is measurable with respect to $\sigma(U,(|X_1^{(j)}|,\dots,|X_d^{(j)}|))$ for $j=1,\dots,m$. Then, for any $\alpha \in (0,0.1]$ there exists $c > 0$ such that
\begin{equation}
\underset{f \in H_\rho}{\sup} \P_f \left( T_{\alpha} = 0 \right) \geq 3/4.
\end{equation}
whenever 
\begin{equation}\label{eq:invariant_rate}
\rho^2 \leq c \frac{\sqrt{d}}{\sqrt{m}n}.
\end{equation}
\end{theorem}
\begin{proof}
In view of Lemma \ref{lem:workhorse} and the proof of the main theorems in \ref{ssec:main_theorem_proof}, it suffices to bound the trace of $\Xi_u$ in \eqref{eq:private_coin_chisq_divergence} and \eqref{eq:public_coin_chisq_divergence} in Lemma \ref{lem:chi-square-bound} (the first term in the exponent is controlled by \eqref{eq:invariant_rate}). 
By assumption on $S^{(j)}$, we have
\begin{equation}\label{eq:sigma_algebra}
\sigma(S^{(j)},U,(|X_1^{(j)}|,\dots,|X_d^{(j)}|) = \sigma(U,(|X_1^{(j)}|,\dots,|X_d^{(j)}|)),
\end{equation}
which implies that the $\text{sign}(X^{(j)}_i)$ is independent of $\sigma(S^{(j)},U,(|X_1^{(j)}|,\dots,|X_d^{(j)}|)$. Writing $X^{(j)}_i = \text{sign}(X^{(j)}_i) |(X^{(j)}_i)|$, we obtain that
\begin{align*}
\E_0\left[ {X}^{(j)} \bigg| S^{(j)}, U=u\right] &= \left( \E_0\left[ \text{sign}(X^{(j)}_i) |(X^{(j)}_i)| \bigg| S^{(j)}, U=u\right]  \right)_{1\leq i \leq d} \\
&= \left( \E_0 \text{sign}(X^{(j)}_i) \E_0\left[ |(X^{(j)}_i)| \bigg| S^{(j)}, U=u\right]  \right)_{1\leq i \leq d} = 0,
\end{align*}
where the second last inequality follows from the fact that $\text{sign}(X^{(j)}_i)$ is independent of the sigma algebra in \eqref{eq:sigma_algebra} and the final equality by the symmetry of the Gaussian distribution around the mean. Following the proof of Theorem \ref{thm:private_coin_lb} with $\Xi_u = 0$, we obtain that the testing risk is bounded from below whenever $\rho^2 \lesssim \frac{\sqrt{d}}{\sqrt{m}n}$.
\end{proof}

\subsection{Lemmas related to rate attainability}\label{ssec:rates_attainable}

\begin{lemma}\label{lem:sqrtm_rate_attained}
Let $T_\alpha$ correspond to a test of level $\alpha$ based on Edginton's method based for p-values $p^{(j)} = \chi^2_d(\| \sqrt{n} X^{(j)}\|_2^2)$ or simply the sum of $\| \sqrt{n} X^{(j)}\|_2^2$. For all $\alpha,\beta \in (0,1)$ if 
\begin{equation}
\rho^2 \geq C_{\alpha,\beta}  \frac{ \sqrt{d}}{\sqrt{m} n}
\end{equation}
we have
\begin{equation*}
\underset{f \in H_\rho}{\sup} \P_f \left( T_\alpha = 0 \right) \leq \beta
\end{equation*}
for $d \geq C_{\alpha,\beta} m$ a large enough constant $C_{\alpha,\beta}$ depending only on $\alpha,\beta \in (0,1)$. The above result holds for Fisher's method also, under the additional assumption that $\log(m) \lesssim \sqrt{d}$.
\end{lemma}
\begin{proof}
The test in \eqref{test:large:d} has level $\alpha$ under the null hypothesis. Under the alternative hypothesis,
\begin{equation*}
\| \sqrt{n} X^{(j)}\|_2^2 \overset{d}{=} n \| f \|_2^2 + 2 \sqrt{n} (Z^{(j)})^\top f + \|Z^{(j)}\|_2^2,
\end{equation*}
where $Z^{(j)} \sim N(0,I_d)$. Rearranging, the test $T_\alpha$ of \eqref{test:large:d} can be seen to equal
\begin{equation}\label{eq:chi_square_alt}
\mathbbm{1}\left\{  2 \frac{\sqrt{n}}{\sqrt{d}} \left(m^{-1/2}\sum_{j=1}^m Z^{(j)}
\right)^\top f + \frac{1}{\sqrt{md}} \underset{j=1}{\overset{m}{\sum}} \left( \| Z^{(j)} \|_2^2 -d\right) \geq \eta_{d,m} - \frac{\sqrt{m}n}{\sqrt{d}} \| f \|_2^2 \right\}
\end{equation}
in distribution under $\P_f$, with
\begin{equation*}
\eta_{d,m} :=\frac{1}{\sqrt{dm}} \left(F_{\chi_{dm}^2}^{-1} (1 - \alpha) - md\right).
\end{equation*}
By Lemma \ref{lem:inverse_CDF_CLT}, $\eta_{d,m} \to \Phi^{-1}(1 - \alpha)$ as both or either $d,m \to \infty$, so $\eta_{d,m}$ is bounded in $d$ and $m$. Consequently, $\P_f(1-T_\alpha)$ equals
\begin{align*}
 \text{Pr}\left( (1 + \frac{\sqrt{n}}{\sqrt{d}} \|f \|_2) O_P(1) \leq \eta_{d,m} - \frac{\sqrt{m}n}{\sqrt{d}} \| f \|_2^2 \right)
\end{align*}
as the left hand side of the test in \eqref{eq:chi_square_alt} is mean $0$ and has constant variance. Since $\|f\|_2^2 \geq C_{\alpha,\beta} \frac{\sqrt{d}}{\sqrt{m}n}$, the latter display can be bounded from above by
\begin{equation*}
\text{Pr}\left( (1 + \frac{\sqrt{n}}{\sqrt{d}} \|f \|_2) O_P(1) \leq - \frac{\sqrt{m}n}{2\sqrt{d}} \| f \|_2^2 \right)
\end{equation*}
for a large enough $C_{\alpha,\beta}$. The latter display is smaller than $\beta$ for $C_{\alpha,\beta} > 0$ large enough depending only on $\alpha$ and $\beta$. 

For Edgington's method, one can take $p^{(j)} = 1 - F_{\chi^2_d}(\| \sqrt{n} X^{(j)}\|_2^2)$ and compute the test
\begin{equation}\label{eq:Edgington_test}
T_\alpha := \mathbbm{1} \left\{ m^{-1/2} \zeta_{\alpha,m} \underset{j=1}{\overset{m}{\sum}} (p^{(j)} - \frac{1}{2})  \geq 12^{-1/2}\Phi^{-1}(1-\alpha) \right\},
\end{equation}
where $\zeta_{\alpha,m} \to 1$ in $m$ is such that $\P_0 T_\alpha = \alpha$, by Lemma \ref{lem:inverse_CDF_CLT}. 

Under the alternative, $\E_f p^{(j)} = \text{Pr}(\| \sqrt{n}f + Z^{(j)}\|_2^2 \leq \chi_{d}^2)$. Therefore, by Lemma 4 in \cite{szabo2022optimal_IEEE}, 
\begin{equation*}
\E_f p^{(j)} \geq \frac{1}{2} + \frac{1}{40} \left( d^{-1/2} n \|f\|_2^2 \bigwedge \frac{1}{2} \right),
\end{equation*}
where we note that we can take $d$ larger than an arbitrary constant as the rate $\sqrt{d}/(\sqrt{m}n)$ being optimal ($\sqrt{d}/(\sqrt{m}n) \lesssim d/({m}n)$) implies $d \gtrsim m$ and for constant order $m$ there is nothing to prove. We obtain that
\begin{align*}
\P_f (1-T_\alpha) &= \P_f \left( \frac{\zeta_{m,\alpha}}{\sqrt{m}}\underset{j=1}{\overset{m}{\sum}} (p^{(j)} - 1/2)  \leq  12^{-1/2}\Phi^{-1}(1-\alpha) \right) \\ 
&= \P_f\left( \frac{\zeta_{m,\alpha}}{\sqrt{m}}\underset{j=1}{\overset{m}{\sum}} [(p^{(j)} - \E_f p^{(j)}) + \E_f p^{(j)} - \frac{1}{2} ] \leq 12^{-1/2}\Phi^{-1}(1-\alpha) \right) \\
&\leq \text{Pr}\left(  O_P(1)  + \frac{\zeta_{m,\alpha} \sqrt{m}}{40} \left( d^{-1/2} n \|f\|_2^2 \bigwedge \frac{1}{2} \right) \leq 12^{-1/2}\Phi^{-1}(1-\alpha) \right),
\end{align*}
where the $O_P(1)$ term in last equality follows from the fact that $\zeta_{m,\alpha}$ is bounded and the central limit theorem (the $p^{(j)}$'s are bounded and independent still under $\P_f$). If the minimum is taken in $1/2$, the result follows for large enough $m$. If the minimum is taken in the first argument,
\begin{equation*}
\frac{\zeta_{m,\alpha} \sqrt{m}}{40} \left( d^{-1/2} n \|f\|_2^2 \bigwedge \frac{1}{2} \right) \geq \frac{C_{\alpha,\beta}\zeta_{m,\alpha}}{40}
\end{equation*}
so for large enough $C_{\alpha,\beta}$, we obtain that $\P_f (1-T_\alpha)\leq \beta$.

For Fisher's method, the test of level $\alpha$ is given by
\begin{equation}\label{eq:Fisher_test}
T_\alpha := \mathbbm{1} \left\{ \underset{j=1}{\overset{m}{\sum}} - 2 \log p^{(j)}  \geq F^{-1}_{\chi^2_{2m}}(1-\alpha) \right\},
\end{equation}
for the p-value $p^{(j)} := 1 - F_{\chi^2_d}(\| \sqrt{n} X^{(j)}\|_2^2)$ (or equivalently Pearson's method for the p-value $F_{\chi^2_d}(\| \sqrt{n} X^{(j)}\|_2^2)$).

For the Type II error bound, assume first that $ n \|f\|_2^2 \geq 20\sqrt{d}$. We have that $\| Z^{(j)}\|_2^2 \geq d - 5\sqrt{d}$ on an event of probability at least $1 - e^{-5}$, via e.g. Theorem 3.1.1 in \cite{vershynin_high-dimensional_2018}. By using a union and a standard Gaussian concentration inequality, the event
\begin{equation}\label{eq:uniform_event_fisher}
\underset{1 \leq j \leq m}{\max}  \left|2 \frac{\sqrt{n}}{\sqrt{d}} f^\top Z^{(j)}\right| \leq \frac{n}{2\sqrt{d}} \|f\|_2^2,
\end{equation}
has mass at least $1- me^{-n \|f\|_2^2/32} \geq 1- me^{-\sqrt{d}/2}$. On the intersection of these two events,
\begin{align*}
F_{\chi^2_d}(\| \sqrt{n}f + Z^{(j)}\|_2^2) &= \text{Pr}\left( \frac{\chi_d^2 - \|Z^{(j)}\|_2^2}{\sqrt{d}} \leq 2 \frac{\sqrt{n}}{\sqrt{d}} f^\top Z^{(j)} + \frac{n}{\sqrt{d}} \|f\|_2^2  \right) \\
&\geq \text{Pr}\left( \frac{\chi_d^2 - d}{\sqrt{d}} \leq \frac{n}{2 \sqrt{d}} \|f\|_2^2  - 5 \right) \\
&\geq \text{Pr}\left( \frac{\chi_d^2 - d}{\sqrt{d}} \leq 5 \right),
\end{align*}
 where the right-hand side tends to $\Phi(5)$ in $d$ by the central limit theorem. As $\Phi(5) > e^{-2}$, we obtain $- \log p^{(j)} \geq 2$. Since $Z^{(1)},\dots,Z^{(m)}$ are independent, by binomial concentration, there are at least $(3/4)m$ indexes $j=1,\dots,m$ such that $\| Z^{(j)}\|_2^2 \geq d - 5\sqrt{d}$ whilst also satisfying \eqref{eq:uniform_event_fisher} with probability $1 - e^{-\tau m} - me^{-\sqrt{d}/2} $ for some constant $\tau > 0$. Using that we can without loss of generality take $m \geq M_{\alpha, \beta}$ for a constant $M_{\alpha,\beta}>0$ (otherwise the separation rate is effectively the same the one for $m=1$) and since we consider $d \gtrsim m$, we obtain that the event the joint event occurs has mass less than $1 - \beta$. Furthermore, on this event, we have $1-T_\alpha = 0$ for $M_{\alpha,\beta} > 0$ large enough, since
\begin{equation*}
\underset{j=1}{\overset{m}{\sum}} - 2 \log p^{(j)} \geq 4m \cdot (3/4)
\end{equation*}
and by the fact that the chi-square quantile tends to $2m + C_\alpha \sqrt{2m}$ for some constant only depending on $\alpha$, which is less than $4m \cdot (3/4) = 3m$ for $m \geq M_{\alpha,\beta}$.

Assume now that $n \| f\|_2^2 \leq 20 \sqrt{d}$. Consider the following claim: for $d$ large enough it holds that
\begin{equation}\label{eq:claim_Fisher}
- 2 \E_f \log \left( 1 - F_{\chi^2_d}(\| \sqrt{n}f + Z^{(j)}\|_2^2) \right) \geq 2 - e^{- \sqrt{d}/4} + \frac{n\|f\|_2^2}{C\sqrt{d}}
\end{equation}
for a fixed constant $C > 0$. If the claim holds,
\begin{equation*}
\P_f (1-T_\alpha) \leq \P_f\left( \frac{\sqrt{m}n\|f\|_2^2}{C\sqrt{d}}  - \sqrt{m} e^{-c \sqrt{d}} + \frac{1}{\sqrt{m}}\underset{j=1}{\overset{m}{\sum}} -2(\log p^{(j)} - \E_f \log p^{(j)}) \leq \eta_{m,\alpha} \right)
\end{equation*}
with $\eta_{m,\alpha} := \frac{1}{\sqrt{2m}}(F^{-1}_{\chi^2_{2m}}(1-\alpha) - 2m)$. Since the method is rate optimal when $m \lesssim d$, the second term of the LHS in the above display may be assumed to be small. For the third term, note that
\begin{align*}
\E_f (-2 \log p^{(j)} )^2 &= 4 \E_f \log ( 1 - F_{\chi^2_d}(\| \sqrt{n}f + Z^{(j)}\|_2^2)) \\
&\leq 4  \log ( 1 - F_{\chi^2_d}(n \| f\|_2^2 + d)), 
\end{align*}
where the last inequality follows from the log-concavity of $x \mapsto 1- F_{\chi_d^2}(x)$ (see e.g. Theorem 3.4 in \cite{finner1997log}). For $n \| f\|_2^2 \leq 20\sqrt{d}$, the latter quanity is uniformly bounded in $n,m$ and $d$. Since the second moment bounds the variance, this implies that 
\begin{equation*}
\frac{1}{\sqrt{m}}\underset{j=1}{\overset{m}{\sum}} -2(\log p^{(j)} - \E_f \log p^{(j)}) = O_P(1)
\end{equation*}
by the independence of $p^{(j)}$ and $p^{(k)}$ for $k\neq j$. Consequently, for some constant $\tau > 0$,
\begin{equation*}
\P_f (1-T_\alpha) \leq \text{Pr}\left( \frac{\sqrt{m}n\|f\|_2^2}{C\sqrt{d}}  - \sqrt{m} e^{-\tau \sqrt{d}} + O_P(1) \leq \eta_{m,\alpha} \right).
\end{equation*}
Since $\eta_{m,\alpha} \to \Phi^{-1}(1-\alpha)$ by Lemma \ref{lem:inverse_CDF_CLT}, the fact that
\begin{equation*}
\frac{\sqrt{m}n\|f\|_2^2}{C\sqrt{d}} \geq C_{\alpha,\beta}/C
\end{equation*}
for large enough $C_{\alpha,\beta}>0$ depending only on $\alpha$ and $\beta$ and the fact that $m \lesssim d$, we have that $\P_f (1-T_\alpha) \leq \beta$.

It remains to prove the claim of \eqref{eq:claim_Fisher}. We start by writing $-2\E_f \log \left( p^{(j)} \right)$ as
\begin{align*}
 -2 \E_f \log (1 - F_{\chi^2_d}(\|Z^{(j)}\|_2^2) -2 \E_f \log \left( \frac{1 - F_{\chi^2_d}(\|\sqrt{n} f + Z^{(j)}\|_2^2)}{1 - F_{\chi^2_d}(\|Z^{(j)}\|_2^2)} \right).
\end{align*}
The first term equals $2$. Using $\log(x) \leq |x - 1|$, the second term is bounded from below by
\begin{equation*}
2 \E_f \left| \frac{F_{\chi^2_d}(\|\sqrt{n} f + Z^{(j)}\|_2^2)-F_{\chi^2_d}(\|Z^{(j)}\|_2^2) }{1 - F_{\chi^2_d}(\|Z^{(j)}\|_2^2)}\right| \geq 2 \E_f \left| F_{\chi^2_d}(\|\sqrt{n} f + Z^{(j)}\|_2^2) - F_{\chi^2_d}(\| Z^{(j)}\|_2^2)\right|. 
\end{equation*}
By the same argument as used for \eqref{eq:uniform_event_fisher},
\begin{equation*}
\E_f \mathbbm{1}_{\{ f^\top Z^{(j)} < 0 \}} F_{\chi^2_d}(\|\sqrt{n} f + Z^{(j)}\|_2^2) \geq \E_f F_{\chi^2_d}(\frac{1}{2}\|\sqrt{n} f\|_2^2 + \|Z^{(j)}\|_2^2)-e^{-\sqrt{d}/4},
\end{equation*}
which is larger than $\E_f F_{\chi^2_d}(\| Z^{(j)}\|_2^2)$ for all large enough $d$. Additionally, on the event that $f^\top Z^{(j)} \geq 0$, it holds that 
\begin{equation*}
 F_{\chi^2_d}(\|\sqrt{n} f + Z^{(j)}\|_2^2) \geq \E_f F_{\chi^2_d}( \frac{1}{2}\|\sqrt{n} f\|_2^2 + \|Z^{(j)}\|_2^2) \geq \E_f F_{\chi^2_d}(\| Z^{(j)}\|_2^2).
\end{equation*}
Furthermore, we have
\begin{equation*}
\E_f F_{\chi^2_d}(\frac{1}{2}\|\sqrt{n} f\|_2^2 + \|Z^{(j)}\|_2^2) - \E_f F_{\chi^2_d}(\| Z^{(j)}\|_2^2) = \text{Pr}\left( 0 \leq \frac{\chi_d^2 - \tilde{\chi}_d^2}{\sqrt{d}} \leq \frac{n}{2\sqrt{d}}\| f\|_2^2 \right),
\end{equation*}
where $\chi_d^2,\tilde{\chi}_d^2$ are independent chi square random variables with $d$ degrees of freedom, which tends in $d$ to 
\begin{equation*}
\Phi\left(\frac{n}{2\sqrt{d}}\| f\|_2^2\right) - \Phi(0) \geq \frac{n}{C\sqrt{d}}\| f\|_2^2,
\end{equation*}
where the inequality holds under the assumption $n \| f\|_2^2 \leq 20 \sqrt{d}$ for a large enough constant $C>0$. Putting the above lower bounds together, we obtain \eqref{eq:claim_Fisher}.
\end{proof}

\begin{lemma}\label{lem:priv_coin_coordinate_strat_rate}
Let $T_\alpha$ correspond to a test of level $\alpha$ considered in \eqref{eq:private_coin_test} or \eqref{eq:edginton_like_test}. For all $\alpha,\beta \in (0,1)$ if
\begin{equation}
\rho^2 \geq C_{\alpha,\beta}  \frac{{d^{3/2}}}{m n}
\end{equation}
we have
\begin{equation*}
\underset{f \in H_\rho}{\sup} \P_f \left( T_\alpha = 0 \right) \leq \beta
\end{equation*}
for a large enough constant $C_{\alpha,\beta}$ depending only on $\alpha,\beta \in (0,1)$.
\end{lemma}
\begin{proof}
The proof follows a similar line of reasoning as e.g. the proof of Lemma A.8 in \cite{szabo2023optimal}. Starting with \eqref{eq:private_coin_test}, note that
\begin{equation*}
\P_f(1-T_\alpha) = \text{Pr} \Big( \frac{1}{\sqrt{d}} \underset{i=1}{\overset{d}{\sum}} \Big(  ( d^{-1/2} \sqrt{mn}f_i + Z_i) \Big)^2  \leq d^{-1/2}F^{-1}_{\chi_d^2}(1-\alpha) \Big)
\end{equation*}
for independent $Z_1,\dots,Z_d \sim N(0,1)$. The latter display equals
\begin{align*}
\text{Pr} \Big( \frac{nm}{d\sqrt{d}} \|f\|_2^2 + 2 \frac{\sqrt{mn}}{d} \underset{i=1}{\overset{d}{\sum}}  f_i Z_i+ \frac{1}{\sqrt{d}} \underset{i=1}{\overset{d}{\sum}}  \left( Z_i^2 - 1 \right) \leq d^{-1/2}(F^{-1}_{\chi_d^2}(1-\alpha) - d) \Big) &= \\
\text{Pr} \Big( \frac{nm}{d\sqrt{d}} \|f\|_2^2 + (1 + \sqrt{\frac{nm}{d^2}} \|f\|_2) O_P(1)  \leq d^{-1/2}(F^{-1}_{\chi_d^2}(1-\alpha) - d) \Big) &\leq \\
\text{Pr} \Big( (1 + \sqrt{\frac{nm}{d^2} \|f\|_2^2}) O_P(1)  \leq - \frac{1}{2} \frac{nm}{d\sqrt{d}} \|f\|_2^2  \Big), &
\end{align*}
where the last inequality holds for large enough $C_{\alpha,\beta}$ since $\frac{nm}{d\sqrt{d}} \|f\|_2^2 \geq C_{\alpha,\beta}$ and $d^{-1/2}(F^{-1}_{\chi_d^2}(1-\alpha) - d)$ is bounded in $d$ by Lemma \ref{lem:inverse_CDF_CLT}. The resulting probability can be made arbitrarily small by taking large enough $C_{\alpha,\beta}$.

For a variation to Edgington's method, ie \eqref{eq:edginton_like_test}, similar reasoning applies. Under the null hypothesis, $\E_0 \Phi (\sqrt{n}X^{(j)}) = 1/2$, so a conservative test (i.e. $\P_0 T_\alpha \leq \alpha$) based on Edgington's method is given by
\begin{equation*}
T_\alpha = \mathbbm{1}\left\{ \left|\frac{1}{\sqrt{d}} \underset{i=1}{\overset{d}{\sum}} \left[ \frac{d}{m}\Big( \underset{j \in \cJ_i}{\overset{}{\sum}}  (p^{(j)} - \frac{1}{2})\Big)^2 - \text{Var}_0(p^{(j)}) \right] \right|  \geq c \alpha^{-1/2} \right\}
\end{equation*}
for a constant $c > 0$ by e.g. Chebyshev's inequality. Under the alternative hypothesis, we have $p^{(j)} = \Phi (\sqrt{n}X^{(j)}_i) = \Phi (\sqrt{n}f_i + Z^{(j)}_i)$ whenever $j \in \cJ_i$. The Type II error $\P_f(1 - T_\alpha)$ equals
\begin{align}
  \P_f\left( \bigg| \frac{1}{\sqrt{d}} \underset{i=1}{\overset{d}{\sum}} \left[ \frac{d}{m}\Big( \underset{j \in \cJ_i}{\overset{}{\sum}}  (p^{(j)} - \Phi (Z^{(j)}_i) + \Phi ( Z^{(j)}_i) - \frac{1}{2})\Big)^2 - \text{Var}_0(p^{(j)}) \right] \bigg|  \leq c\alpha^{-1/2} \right)&= \nonumber \\
 \P_f\left( \bigg| \zeta + \xi + \frac{1}{\sqrt{d}} \underset{i=1}{\overset{d}{\sum}} \frac{d}{m}\Big( \underset{j \in \cJ_i}{\overset{}{\sum}}  (\Phi (\sqrt{n}f_i + Z^{(j)}_i) - \Phi(Z^{(j)}_i)) \Big)^2 \bigg|  \leq c\alpha^{-1/2} \right) \label{eq:high_dim_strat_continue_here}
\end{align}
where 
\begin{equation*}
\zeta = \frac{1}{\sqrt{d}} \underset{i=1}{\overset{d}{\sum}} \left[ \frac{d}{m}\Big( \underset{j \in \cJ_i}{\overset{}{\sum}}  (\Phi (Z^{(j)}_i) - \frac{1}{2}) \Big)^2 - \text{Var}_0(p^{(j)}) \right]
\end{equation*}
and
\begin{equation*}
\xi = \frac{2\sqrt{d}}{m} \underset{i=1}{\overset{d}{\sum}} \Big( \underset{j \in \cJ_i}{\overset{}{\sum}}  (\Phi (\sqrt{n}f_i + Z^{(j)}_i) - \Phi (Z^{(j)}_i)) \Big) \Big( \underset{j \in \cJ_i}{\overset{}{\sum}} (\Phi (Z^{(j)}_i) - \frac{1}{2}) \Big).
\end{equation*}
By independence between $Z^{(j)}_i$ and $Z^{(k)}_i$ when $j \neq k$, the random variable $\zeta$ is mean $0$ under $\E_f$ with constant variance (i.e. not depending on $d,m,n$) and is thus $O_P(1)$. Similarly, $\xi$ has constant order variance and expectation. 
By Jensen's inequality
\begin{equation*}
 \E_f (\Phi (\sqrt{n}f_i + Z^{(j)}_i) - \Phi(Z^{(j)}_i))^2 \geq (\Phi (2^{-1/2}\sqrt{n}f_i) - \Phi(0))^2
\end{equation*}
where it is used that 
\begin{equation*}
\E_f \Phi (\sqrt{n}f_i + Z^{(j)}_i) = \text{Pr} \left( \sqrt{n}f_i + Z \geq Z' \right) = \Phi( 2^{-1/2}\sqrt{n}f_i ).
\end{equation*}
By Lemma A.11 in \cite{szabo2023optimal}, the RHS of the second last display is lower bounded by $\frac{1}{12} \min \{ \frac{1}{2}n f_i^2, 1 \}$. By the independence of $Z^{(j)}_i$ and $Z^{(k)}_i$ when $j\neq k$, it also holds that
\begin{align*}
\E_f(\Phi (\sqrt{n}f_i + Z^{(j)}_i) - \Phi(Z^{(j)}_i)) (\Phi (\sqrt{n}f_i + Z^{(k)}_i) - \Phi(Z^{(k)}_i)) =(\Phi (2^{-1/2}\sqrt{n}f_i) - \Phi(0))^2.
\end{align*}
Therefore,
\begin{equation*}
\E_f \frac{d}{m}\Big( \underset{j \in \cJ_i}{\overset{}{\sum}}  (\Phi (\sqrt{n}f_i + Z^{(j)}_i) - \Phi(Z^{(j)}_i)) \Big)^2 \geq \frac{m}{12d} \min \{ \frac{1}{2}n f_i^2, 1 \}.
\end{equation*}
Adding and subtracting the above expectation and noting that 
\begin{equation*}
\frac{1}{\sqrt{d}} \underset{i=1}{\overset{d}{\sum}} \frac{d}{m}\Big( \underset{j \in \cJ_i}{\overset{}{\sum}}  (\Phi (\sqrt{n}f_i + Z^{(j)}_i) - \Phi(Z^{(j)}_i)) \Big)^2
\end{equation*}
has constant variance by the independence of $Z^{(j)}_i$ and $Z^{(k)}_i$ when $j\neq k$, we obtain that \eqref{eq:high_dim_strat_continue_here} is bounded above by
\begin{equation*}
\P_f\left(  O_P(1) + \frac{m}{12 d \sqrt{d}} \underset{i=1}{\overset{d}{\sum}}  \min \{ \frac{1}{2}n f_i^2, 1 \}   \leq c\alpha^{-1/2} \right).
\end{equation*}
If the minimum is taken by $1$ for any $i=1,\dots,d$, the proof is completed by noting that $m \gtrsim d^2$ by assumption whenever the rate $\frac{d \sqrt{d}}{nm}$ is the optimal rate and considering $m$ large enough. Otherwise, the power is arbitrarily small for 
\begin{equation*}
\frac{mn}{ d \sqrt{d}} \|f\|_2^2 \geq C_{\alpha,\beta}
\end{equation*}
and $C_{\alpha,\beta}$ large enough.
\end{proof}

\begin{lemma}\label{lem:pub_coin_coordinate_strat_rate}
Let $T_\alpha$ correspond to a test of level $\alpha$ considered in \eqref{eq:pub_coin_test}. For all $\alpha,\beta \in (0,1)$ if
\begin{equation}
\rho^2 \geq C_{\alpha,\beta}  \frac{{d}}{m n}
\end{equation}
we have
\begin{equation*}
\underset{f \in H_\rho}{\sup} \P_f \left( T_\alpha = 0 \right) \leq \beta
\end{equation*}
for a large enough constant $C_{\alpha,\beta}$ depending only on $\alpha,\beta \in (0,1)$.
\end{lemma}
\begin{proof}
The proof follows a similar line of reasoning as e.g. the proof of Lemma A.7 in \cite{szabo2023optimal}. For any $f \in \R^d$ such that $\|f\|_2 \geq \rho$, we have
\begin{equation*}
U \sqrt{n}X^{(j)} \overset{d}{=} \sqrt{n} U f + Z^{(j)}
\end{equation*}
under $\P_f$ by rotational invariance of the normal distribution. The probability of a Type II error of the test of level $\alpha$ given in \eqref{eq:pub_coin_test} is then equal to
\begin{equation*}
\text{Pr}\left( \big| \sqrt{n} \sqrt{m} (U f)_1 + Z \big| \leq \Phi^{-1}(1-\alpha/2)  \right),
\end{equation*}
with $Z \sim N(0,1)$. The random variable $(U f)_1$ is in distribution equal to $\|f\|_2 Z_1' / \|Z'\|_2$ for a $d$-dimensional standard Gaussian random vector $Z'$. For any $\beta \in (0,1)$, there exists $c' > 0$ such that $\|Z'\|_2 > c' \sqrt{d}$ occurs with probability $1 - \beta/2$. Also, for $\frac{\sqrt{nm} \|f\|_2}{c' \sqrt{d}} \geq C_{\alpha,\beta}/c'$ large enough, 
\begin{equation*}
\text{Pr}\left( \big| \frac{\sqrt{nm} \|f\|_2}{c' \sqrt{d}}  + Z \big| \leq \Phi^{-1}(1-\alpha/2)  \right) \leq \beta/2.
\end{equation*}
This concludes the proof of the lemma.
\end{proof}

The following fact is well known and included for completeness. For a random variable $V$, let $F_V$ denote its CDF. 
\begin{lemma}\label{lem:inverse_CDF_CLT}
Let $W_1,\dots,W_m$ be random variables and let $V_m = \sum_{j=1}^{m} W_j$. Suppose that 
\begin{equation*}
m^{-1/2} \underset{j=1}{\overset{m}{\sum}} (W_j - \E W_j) \rightsquigarrow N(0, \sigma^2).
\end{equation*}
Then, for all $\alpha \in (0,1)$,
\begin{equation*}
(\sigma^2 m)^{-1/2} \left(F_{V_m}^{-1} (\alpha) - \sum_{j=1}^{m} \E W_j \right) \to \Phi^{-1}(\alpha),
\end{equation*}
where $\Phi$ is the standard Gaussian CDF.
\end{lemma}
\begin{proof}
The quantile function 
\begin{equation*}
F_{V_m}^{-1} (\alpha) = \inf \left\{ x \in \R : \text{Pr}\left( V_m \leq x \right) \geq \alpha \right\}
\end{equation*}
satisfies $z(F_{V_m}^{-1} (\alpha) - y) = F_{z(V_m - y)}^{-1} (\alpha)$. The result now follows by e.g. Lemma 21.2 in \cite{vaart_asymptotic_1998}.
\end{proof}

\subsection{Proof Lemma \ref{lem:p_value_combs} and Lemma \ref{lem:e_value_combs}}\label{ssec:proofs_p_e_values}

\begin{proof}[Proof of Lemma \ref{lem:p_value_combs}]
The lemma directly follows from Theorem \ref{thm:private_coin_lb} and Theorem \ref{thm:public_coin_lb} after verifying the corresponding conditions. Assumption \ref{assumption:local_randomness_statistic} is satisfied if $p^{(j)}$ is generated using only local randomness, while in case of shared randomness, the same conclusion holds for Assumption \ref{assumption:public_randomness_statistic}. Below, we prove Assumptions  \ref{assumption:lipschitz_type_C_m} and \ref{assumption:test_form} for the examples listed in the lemmas.

\begin{enumerate}
    \item  Fisher's method: let $S^{(j)} = -2 \log p^{(j)}\sim^{H_0}\chi_2^2$ and consider the test of level $\alpha$ as
    \begin{equation*}
    \mathbbm{1}\Big \{ \eta_{\alpha,m} \frac{1}{\sqrt{2m}} \underset{j=1}{\overset{m}{\sum}} (S^{(j)} - 2)  \geq \Phi^{-1}({1-\alpha}) \Big\}
    \end{equation*}
    with $\Phi^{-1}$ the inverse standard normal CDF and
    \begin{equation*}
        \eta_{\alpha,m} := \Phi^{-1}({1-\alpha}) \left(\frac{1}{\sqrt{2m}}\Big(F_{\chi^2_{2m}}^{-1}(1-\alpha) - 2m\Big) \right)^{-1}. 
        \end{equation*}
In view of the CLT, see Lemma \ref{lem:inverse_CDF_CLT}, the sequence $\eta_{\alpha,m}$  converges to one, hence it is bounded. Furthermore, note that the corresponding combination function $C_m(s) := (\eta_{\alpha,m}/\sqrt{m}) {\sum}_{j=1}^m (s_j - 1)$ with  $s=(s_j) \in \R^m$ satisfies Assumption \ref{assumption:lipschitz_type_C_m} (e.g. with $p=q=1$). This in turn implies the moment condition for $S^{(j)}$, concluding the proof.
    \item Mudholkar and George's method: The corresponding combination function $C_m(s) := |m^{-1/2}\sum_{j=1}^m s_j |$, by triangle inequality, satisfies Assumption \ref{assumption:public_randomness_statistic}. Since $S^{(j)} := -\log (p^{(j)}(1-p^{(j)}))$, the moment conditions are also satisfied. 
    \item Pearson's and Edgington's methods: the proofs follow the same reasoning as above with an additional application of the reverse triangle inequality in case of a two sided test.
    \item Tippett's method: when small p-values are expected under the alternative hypothesis, a test of level $\alpha \in (0,1)$ is given by
    \begin{equation*}
    T_\alpha = \mathbbm{1}\left\{ 1 - \big(1-\min \{ p^{(1)},\dots,p^{(m)} \}\big)^m \leq \alpha \right\}, 
    \end{equation*}
    where $1 - (1-\min \{ p^{(1)},\dots,p^{(m)} \})^m$ is uniformly distributed under the null (see e.g. \cite{tippett1941methods}). Observe that it is equivalent to
    \begin{equation*}
    \mathbbm{1}\left\{  - m  \min \left\{ - \log (1-p^{(j)})\right\} \geq  \log (1-\alpha) \right\}.
    \end{equation*}
    For $j=1,\dots,m$, take $S^{(j)} = - \log (1-p^{(j)})\sim^{H_0}\text{Exp}(1)$. The threshold $\alpha \mapsto \log (1-\alpha)$ is strictly decreasing and the combination function $C_m(s) = - m \min s_j$ satisfies
    \begin{equation*}
    |C_m(s) - C_m(s')| \leq m \min |s_j - s_j'| \leq \underset{j=1}{\overset{m}{\sum}} |s_j - s_j'|. 
    \end{equation*}
    Consequently, Assumptions \ref{assumption:test_form} and \ref{assumption:lipschitz_type_C_m} are satisfied.
    \item Generalized averages: The case where $r = -\infty$ corresponds to Tippett's method above. Similarly, $r = \infty$ corresponds to the maximum of p-values, for which the proof follows by similar steps. For $r \in [\frac{1}{m-1}, \infty)$, $a_{r,m}$ can be chosen such that the test $ T_\alpha$ in defined in Section \ref{ssec:p-values_combined} has precise level: $\P_0 T_\alpha = \alpha$, see Proposition 2 and 3 in \cite{vovk2020combining}. For such $a_{r,m}$, the set $\{ a_{r,m} : r \in [\frac{1}{m-1}, \infty) \}$ is bounded (see Table 1 in the aforementioned paper). This test can easily be seen to be of the form \eqref{assumption:test_form} and for the generalized average, we have
    $ ( m^{-1} \sum_{j=1}^m (s_j)^r  )^{1/r} = \| m^{-1/r}s \|_r$, which yields
\begin{align*}
 m^{-1/r} \big|  \| s \|_r - \| s' \|_r  \big| \leq  m^{-1/r} \| s - s' \|_r 
 \leq \max_j |s_j - s_j'|,
\end{align*}
so Assumption \ref{assumption:lipschitz_type_C_m} is satisfied since $a_{r,m}$ is bounded. 
\end{enumerate}

\end{proof}

\begin{proof}[Proof of Lemma \ref{lem:e_value_combs}]
\textbf{Product of e-values:} The e-value test $T_\alpha$ for the combination function $(e_j) \mapsto \underset{j=1}{\overset{m}{\Pi}} e_j $ can be written as
\begin{equation*}
T_\alpha = \mathbbm{1}\Big\{ \underset{j=1}{\overset{m}{\sum}}  \log E^{(j)} \geq \log(1/\alpha) \Big\}.
\end{equation*}
 For $S^{(j)} := \log E^{(j)}$ and $C_m(s) =\sum_{j=1}^m s_j$ note that $E_0 | \log E^{(j)}| < \infty$ and $C_m$ satisfies \eqref{eq:C_m_approximation_assumption}. Since $\alpha \mapsto \log(1/\alpha)$ is strictly decreasing on $(0,1)$, the assumptions of Theorems  \ref{thm:private_coin_lb} and \ref{thm:public_coin_lb} are met.

\textbf{Average of e-values:} Since $E^{(j)}$ is nonnegative, the moment condition is satisfied. The map $ (e_j) \mapsto m^{-1} \sum_{j=1}^m e_j $ satisfies \eqref{eq:C_m_approximation_assumption}, while the map $\alpha \mapsto \alpha^{-1}$ is strictly decreasing and independent of $m$. Hence the conditions of  Theorems  \ref{thm:private_coin_lb} and \ref{thm:public_coin_lb} are satisfied.
\end{proof}

\subsection{Additional simulations}\label{sec:additional_sims}

\begin{figure}[ht]\label{fig:roc_curves_d}
\includegraphics[width=0.9\textwidth]{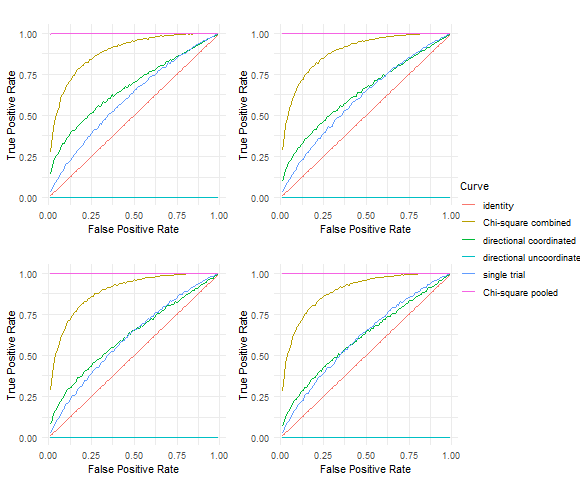}
\caption{ROC curves for different values of $d$, whilst keeping $m=20$, $n=30$, $\rho^2 =  9\sqrt{d}/(16n)$. From left to right, top to bottom: $d=30$, $d=60$, $d=90$, $d=120$. The uncoordinated directional test requires $m \geq d$ and is therefore has TPR set to $0$.}
\end{figure}

\begin{figure}[ht]\label{fig:roc_curves_m}
\includegraphics[width=0.9\textwidth]{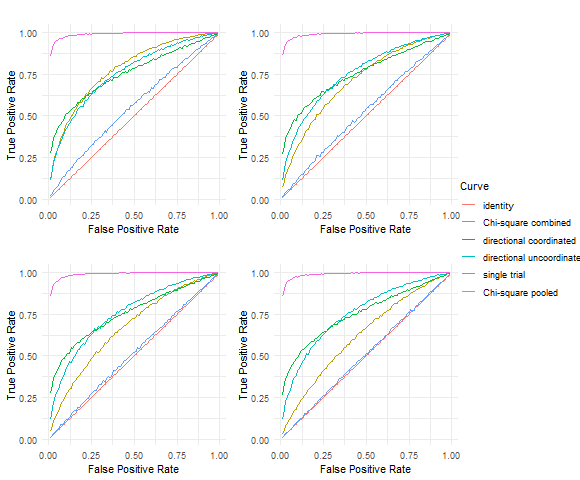}
\caption{ROC curves for different values of $m$, whilst keeping $d=5$, $n=30$, $\rho^2 =  9d/(16nm)$. From left to right, top to bottom: $m=30$, $m=60$, $m=100$, $m=200$.}
\end{figure}

Figure \ref{fig:roc_curves_d} shows the further improvement of the combined chi-square tests compared to the directional methods as $d$ grows with respect to the number of trials, for signals that are around the detection threshold. Figure \ref{fig:roc_curves_m} shows the further worsening of performance of the combined chi-square tests compared method as $m$ grows with respect to the dimension, for signals that are around the detection threshold. For each of these simulations, $10,000$ repitions for every value $\alpha \in \{ 0.01, 0.02,\dots, 0.99 \}$ of the level of the tests are considered.

\end{document}